\documentclass[10pt]{article}
\usepackage{amsmath, amssymb,epsfig,bm}
\usepackage{latexsym}
\usepackage{mathtools}
\usepackage{graphicx}
\usepackage{color}
\numberwithin{equation}{section}
\DeclareMathAlphabet{\itbf}{OML}{cmm}{b}{it}

\newcommand{\RR}{\mathbb{R}}

\newcommand{\ds}{\displaystyle}

\def\debproof{ {\bf Proof.} }
\def\finproof{\hfill {\small $\Box$} \\}

\newcommand{\bp}{{\itbf p}}
\newcommand{\bq}{{\itbf q}}
\newcommand{\br}{{\itbf r}}
\newcommand{\bs}{{\itbf s}}
\newcommand{\bz}{{\itbf z}}
\newcommand{\bx}{{\itbf x}}
\newcommand{\bv}{{\itbf v}}

\newcommand{\bu}{{\itbf u}}
\newcommand{\bn}{{\itbf n}}
\newcommand{\by}{{\itbf y}}

\newcommand{\bR}{{\itbf R}}
\renewcommand{\i}{\mathrm{i}}

\newcommand{\bE}{{\itbf E}}
\newcommand{\bH}{{\itbf H}}
\newcommand{\bJ}{{\itbf J}}
\newcommand{\bF}{{\itbf F}}
\newcommand{\bG}{\itbf G}

\newcommand{\be}{\begin{eqnarray}}
\newcommand{\ee}{\end{eqnarray}}
\newcommand{\nn}{\nonumber}
\newcommand{\ben}{\begin{eqnarray*}}
\newcommand{\een}{\end{eqnarray*}}
\def\ds{\displaystyle}
\def\nm{\noalign{\medskip}}

\newtheorem{lem}{Lemma}[section]
\newtheorem{prop}{Proposition}[section]
\newtheorem{cor}{Corollary}[section]
\newtheorem{thm}{Theorem}[section]
\newtheorem{rem}{Remark}[section]

\newtheorem{definition}{Definition}[section]
\newtheorem{alg}{Algorithm}[section]

\newcommand{\rot}{{\nabla\times}}

\newcommand{\ddiv}{{\nabla\cdot}}

\begin{document}

\title{Detection and classification from electromagnetic induction data\thanks{\footnotesize This work was supported  by
ERC Advanced Grant Project MULTIMOD--267184, China NSF under the
grants 11001150,  41230210, and 11021101, and National Basic
Research Project under the grant 2011CB309700.}}

\author{Habib Ammari\thanks{\footnotesize Department of Mathematics and Applications,
Ecole Normale Sup\'erieure, 45 Rue d'Ulm, 75005 Paris, France
(habib.ammari@ens.fr,
han.wang@ens.fr).} \and Junqing Chen\thanks{\footnotesize
Department of Mathematical Sciences, Tsinghua University, Beijing
100084, China (jqchen@math.tsinghua.edu.cn).} \and Zhiming
Chen\thanks{\footnotesize LSEC, Institute of Computational
Mathematics Chinese Academy of Sciences, Beijing 100190, China
(zmchen@lsec.cc.ac.cn).}  \and Darko
Volkov\thanks{Department of Mathematical Sciences, Stratton Hall,
100 Institute Road, Worcester, MA 01609-2280, USA
(darko@wpi.edu).}\and Han Wang\footnotemark[2]}

\maketitle

\begin{abstract}
In this paper we introduce an efficient algorithm for identifying conductive objects using induction data derived from eddy currents. Our method consists of  first extracting geometric features from the induction data and then matching them to precomputed  data for known objects from a given dictionary. The matching step relies on fundamental properties of conductive polarization tensors and new invariants introduced in this paper. A new shape identification scheme is introduced and studied. We test it numerically in the presence of measurement noise.  Stability and resolution 
capabilities
 of the proposed identification algorithm are  quantified in numerical simulations.
\end{abstract}

\bigskip

\noindent {\footnotesize Mathematics Subject Classification
(MSC2000): 35R30, 35B30}

\noindent {\footnotesize Keywords: eddy current imaging, induction data, classification, recognition,
invariant shape descriptors}

\section{Introduction}

Electromagnetic induction sensors operate by emitting magnetic fields and detecting the response from  electric currents generated when these fields interact with metallic objects (often referred to as targets). These sensors comprise a transmission coil and a receiver coil. Electric currents flowing from the transmitter coil radiate to produce a primary magnetic field that penetrates the surrounding medium and any nearby metallic objects. A time-variable primary magnetic field induces so-called eddy currents in surrounding metallic objects, and these currents in turn yield a secondary magnetic field which is then sensed by the receiver coil \cite{induction1,induction2}.

Electromagnetic induction sensors are quite sensitive and can detect buried  land mines of
low metallic content
or unexploded ordnance containing only a few grams of metal. At present, commercially available sensors have a limited  ability to distinguish land mines and unexploded ordnance from  metallic clutter. False alarms generated by metallic clutter severely limit the speed and efficiency of land mine clearance operations \cite{review, torquil}.

So far little is known about how the signals collected by these sensors from land mines and unexploded ordnance depend on  operating frequency and on  shape,  location,  size, and  orientation of  metallic targets \cite{gao, torquil, khadr}. Electromagnetic induction has become, however, the technology of choice for detecting and classifying concealed weapons \cite{torquil}. Most weapons typically contain some amount of metal. Each particular weapon has a characteristic electromagnetic signature determined by its size, shape and material composition. 
Currently, the use of electromagnetics based safety systems in airports, railway stations, courts, and so on, is widespread.
Metal detectors commonly used by security agents are, however, plagued by high false alarm rates. 
This is mainly because they are designed to simply be set off once a threshold for  quantity of metal is reached. This makes it at times difficult to differentiate weapons from everyday items. Additionally,  human bodies can alter the sensitivity of  detectors  since they are
themselves slightly conductive. This can lead to poor reliability of  detection systems and may even cause metallic objects to go undetected.

The aim of this paper is to contribute to technologies based on  electromagnetic induction sensors.
In particular we aim at improving detection, characterization, and classification methods. We propose efficient algorithms to better differentiate between land mines, unexploded ordnance or weapons from harmless metallic objects. We believe that our new methods will lead to a drastic reduction in false alarm rates. Our proposed algorithms are able to quickly, accurately, and robustly detect and classify metallic objects using readings of electromagnetic induction measurements. The electromagnetic object classification problem is by nature very challenging since the dependence of  electromagnetic induction data on shape, location, and orientation of targets is highly nonlinear. 
An additional hurdle is that induction data and other distinguishable geometric features of the objects to be imaged
depend on frequency.

In previous work, \cite{ACCGV}, we introduced a novel mathematical analysis and we presented numerical methods pertaining to
 imaging of arbitrary shaped small-volume conductive objects using electromagnetic induction data. We derived in that paper a small-volume expansion of the eddy current data measured at some distance away from the conductive object. That expansion involves two polarization tensors: one associated to magnetic contrast and another to conductivity. These tensors depend intrinsically on the geometry of objects to be imaged. A subspace projection algorithm was designed for locating spherical objects from multistatic response matrix data at a single frequency. That algorithm is of MUSIC type (MUSIC stands for
MUltiple Signal Classification). It uses projections of magnetic dipoles located at search points onto the image space of multistatic response matrices. The $nm$-th entry of these multistatic response matrix is the signal recorded by the $n$-th receiver as the $m$-th source is emitting.  Multistatic measurements
 were shown to significantly increase detection rates and reduce false alarm rates in the presence of measurement noise \cite{AGJKLSW, AGJ, SIAP}.
In this paper, we first show that conductive polarization tensors can be robustly extracted from induction data. We then  derive 
important scaling, rotation, and translation properties of these conductive polarization tensors. Based on these properties, we construct shape descriptors from multifrequency induction data and we then search for a match within a dictionary of targets. Finally, we numerically quantify the stability of the proposed identification algorithm. Interestingly, we also found out that there are objects that could not be unambiguously identified using single frequency data but that became possible to recognize through the use of multiple frequency data. Our proposed identification algorithm involves two steps. First, the metallic object is detected and its location approximately determined using a subspace location algorithm; second, the conductive polarization tensors at multiple frequencies are extracted from the induction data and shape descriptors. These descriptors are invariant with respect to translation and rotation. After reconstructing them, the shape of the object to be imaged is matched to a shape from our pre computed dictionary. We expect our identification algorithm to outperform any method currently employed  to find land mines, unexploded ordnance. Classification algorithms have been recently introduced in electrolocation \cite{dic1,dic2, ABGKW, ammari_generalized_2011}  and in echolocation \cite{echolocation}. 

This paper is organized as follows. In section \ref{sec2} we summarize the main findings from our previous study
on small volume asymptotic theory for eddy currents. 
\color{black} As many concepts and objects related to  eddy currents in unbounded 
domains were introduced in section \ref{sec2}, in section \ref{new_section}
 we are able to state in a concise fashion what precisely is the detection and classification problem that we propose to solve in this paper.
\color{black}
We then present in section \ref{Localization algorithm} a target subspace localization algorithm. 
Section \ref{Properties of the CPTs} is devoted to scaling, rotation, and translation properties of conductive polarization tensors. In section \ref{CPTs recovery}, we show how to recover conductive polarization tensors from electromagnetic data using a least squares 
minimization method and we introduce a classification algorithm. 
In section \ref{Numerical examples}, we show 
a numerical example of localization disambiguation of targets using our algorithm.
  In the last section we close this paper by giving a few concluding remarks, and pointing to directions for future work.

\section{Asymptotic formula for eddy current equations}\label{sec2}

In this section, we recall the asymptotic formula for the eddy current problem with small-volume target. Such a formula extends  the small-volume framework \cite{lnm, book, khelifi, AVV, AV, vog, kwon_real_2002,VV} for imaging conductive targets. 

 Suppose that there is an electromagnetic target in $\RR^3$ of
the form $B_\alpha=\bz+\alpha B$, where $B \subset \RR^3$ is a bounded, smooth
 domain containing the origin. Let $\Gamma$ and $\Gamma_\alpha$ denote the boundary of $B$
  and $B_\alpha$. Let $\mu_0$ denote the magnetic permeability
of  the free space. Let $\mu_*$ and $\sigma_*$
denote the permeability  and the conductivity of the target which are also assumed to be
constant. We introduce the piecewise constant magnetic
permeability and electric conductivity
\begin{eqnarray*}
\mu_{\alpha}(\bx) = \left \{
    \begin{array}{ll}
    \mu_*  & \mbox  {in $B_\alpha$,   } \\
\nm
        \mu_{0}   & \mbox  {in $B^c_\alpha :=\RR^3 \backslash  \overline{B}_\alpha,$}
           \end{array}
 \right.\ \ \ \
 \sigma_\alpha (\bx)=  \left \{
    \begin{array}{ll}
    \sigma_*  & \mbox  {in $B_\alpha$,  }\\
\nm
        0   & \mbox  {in $B_\alpha^c$.}
           \end{array}
 \right.
\end{eqnarray*}
Let $(\bE_\alpha, \bH_\alpha)$ denote the eddy current fields in
the presence of the electromagnetic target $B_\alpha$ and a
source current $\bJ_0$ located outside the target. Moreover, we
suppose that $\bJ_0$ has a compact support and is divergence free:
$\nabla \cdot \bJ_0 =0$ in $\RR^3$. The fields $\bE_\alpha$ and
$\bH_\alpha$ are the solutions of the following eddy current
equations: \be \label{alpha} \left\{
\begin{array}{l}
\rot \bE_\alpha =\i\omega \mu_\alpha \bH_\alpha \ \ \
\hbox{ in } \RR^3,
\\
\nm
 \rot  \bH_\alpha =~\sigma_\alpha \bE_\alpha + \bJ_0 \ \ \
 \hbox{ in } \RR^3,\\
 \nm
 \bE_\alpha(\bx)=O(|\bx|^{-1}), \ \ \bH_\alpha(\bx)=O(|\bx|^{-1}) \ \ \ \mbox{ as } |\bx|\rightarrow \infty.
\end{array} \right.
\ee By eliminating $\bH_\alpha$ in \eqref{alpha} we obtain the
following $\bE$-formulation of the eddy current problem
\eqref{alpha}: \be \label{eform} \left\{
\begin{array}{l}
\rot\mu^{-1}_\alpha\rot\bE_\alpha-
\i\omega\sigma_\alpha\bE_\alpha= \i\omega\bJ_0\ \ \  \mbox{ in }
\RR^3, \\ \nm \ddiv\bE_\alpha =0\ \ \ \mbox{ in } B_\alpha^c, \\
\nm \bE_\alpha(\bx)=O(|\bx|^{-1})\ \ \ \mbox{ as }
|\bx|\rightarrow\infty.
\end{array}\right.
\ee

We denote by  $\bE_0$ the  \color{black}solution of the problem
\begin{equation} \label{E0}
\left\{
\begin{array}{l}
\rot\mu^{-1}_0\rot\bE_0=\i\omega\bJ_0\ \ \  \mbox{ in } \RR^3 , \\
\nm \ddiv \bE_0 = 0 \mbox{ in } \RR^3 , \ \ \  \\ \nm
\bE_0(\bx)=O(|\bx|^{-1})\ \ \  \mbox{ as } |\bx| \rightarrow
\infty.
\end{array}\right.
\end{equation}
\color{black} Problem (\ref{eform}) has a unique solution in appropriate
functional spaces provided we require the additional condition
$\int_{\Gamma_\alpha} \bE_\alpha^+ \cdot \bn =0$ where $\bn$
is the exterior normal vector on $\Gamma_\alpha$ and 
$\bE_\alpha^+$ is the exterior trace of $\bE_\alpha$
on $\Gamma_\alpha$: we refer the reader to \cite{ACCGV} 
for an in depth study of questions regarding well posedness of 
such eddy current equations in unbounded domains. Note that problem 
(\ref{E0}) can be thought of being the unperturbed case ($\alpha=0$) of problem
(\ref{eform}).\\
For problem (\ref{E0}) we require the additional condition
$\int_{\Gamma} \bE_0^+ \cdot \bn =0$, and we set
$\bH_0=\frac{1}{\i\omega\mu_0}\rot\bE_0$.
\color{black}

Let
$k=\omega\mu_0\sigma_*$. We are interested in the asymptotic regime
when $\alpha \rightarrow 0$ and
\begin{equation} \label{deff5}
\nu := k \alpha^2
\end{equation}
is of order one. Moreover, we assume that $\mu_*$ and $\mu_0$ are
of the same order. In eddy current imaging the wave equation is converted into the
diffusion equation, where the characteristic length is the skin
depth $\delta$, given by $\delta = \sqrt{2/k}$. Hence, in the
regime $\nu =O(1)$, the skin depth $\delta$ 
has same order of magnitude as 
 the
characteristic size $\alpha$ of the target.

We denote by $C$ a
 generic constant which depends possibly on $\mu_*/\mu_0$,  the upper bound of
$\omega \mu_0\sigma_* \alpha^2$,
 the domain $B$, but is otherwise independent of
 $\omega,\sigma_*,\mu_0,\mu_*$.

Let $G(\bx,\by)=\frac 1{4\pi |\bx-\by|}$ be the fundamental
solution of the Laplace equation. Let ${\bm \theta}_i$ be the solution of the following interface problem:
\begin{equation}
\left\{
\begin{array}{l}
\nabla_{\bm \xi} \times \mu^{-1} \nabla_{\bm \xi} \times{\bm \theta}_i-\i \omega\sigma\alpha^2 {\bm
\theta}_i= \displaystyle\i\omega\sigma\alpha^2{\itbf
e}_i\times\bm\xi\mbox{ in }B\cup B^c,\\\nm
 \nabla_{\bm \xi} \cdot {\bm
\theta}_i = 0 \mbox{ in } B^c,\\\nm \mbox{[}{\bm \theta}_i
\times\bn]_{\Gamma} =0,\ \ [\mu^{-1} \nabla_{\bm \xi} \times{\bm \theta}_i\times\bn
]_{\Gamma}= - 2 [\mu^{-1}]_\Gamma {\itbf e}_i\times \bn \mbox{ on } {\Gamma}, \\
\nm {\bm
\theta}_i({\bm \xi})=O(|\bm\xi|^{-1})\mbox{ as } |{\bm \xi}|
\rightarrow\infty,
\end{array}
\right. \label{theta_i_def}
\end{equation}
where $\mu(\bm \xi) = \mu_*$ if $\bm \xi \in B$, $\mu(\bm \xi) = \mu_0$ if $\bm \xi \in B^c$
and $\sigma(\bm \xi) = \sigma_*$ if $\bm \xi \in B$, $\sigma(\bm \xi) = 0$ if $\bm \xi \in B^c$, and let $\itbf{e}_i$ be the unit vector in the $x_i$ direction.
This interface problem is uniquely solvable
\color{black} if we require the additional condition
\begin{equation}\int_{\Gamma} {\bm
\theta}_i^+ \cdot \bn =0; \label{theta_i_unique}
\end{equation}
\color{black} see \cite{ACCGV}.

In \cite{ACCGV}, we have proved the following asymptotic formula.

\begin{thm}\label{thm-1} Assume that $\nu$ is of order one and let $\alpha$ be small. For $\bx$ away from the location $\bz$ of the target, we have
\ben
\bH_\alpha(\bx)-\bH_0(\bx)&=&\i\nu\alpha^3\left[\frac 12\,
\sum^3_{i=1}\bH_0(\bz)_i\int_B{\bm D}^2_\bx G(\bx,\bz)\bm\xi\times(\bm\theta_i
+{\itbf e}_i\times\bm\xi)d\bm\xi\right]\nn\\
\nm &&+\alpha^3\Big(1-\frac{\mu_0}{\mu_*}\Big)
\left[\sum^3_{i=1}\bH_0(\bz)_i {\bm D}^2_\bx G(\bx,\bz)
\int_B\Big({\itbf e}_i+\frac 12\nabla\times\bm\theta_i\Big)
d\bm\xi\right]+\bR(\bx), \een where $({\bm D}^2_\bx G )_{ij}=
\partial_{x_i x_j}^2 G$ and  \ben |\bR(\bx)|\le
C \alpha^4\|\bH_0\|_{W^{2,\infty}(B_\alpha)}, \een uniformly in
$\bx$ in any compact set away from $\bz$.

\end{thm}

\color{black}
\begin{definition}
\color{black}
For an arbitrary shaped target $B$ with conductivity $\sigma$ and
size $\alpha$, 
and for $l,l^\prime=1,2,3$,
we define the conductivity
polarization tensor (CPT) $\mathbb{M}^{l,l^\prime}$ 
  to  be the
$3\times 3$ matrix whose $i$-th column is
\be\ds 
\mathbb{M}^{l,l^\prime}_i =
\frac 12
 {\itbf e}_l   \times \int_B \xi_{l^\prime} (\bm\theta_i+{\itbf
e}_i\times\bm\xi)d\bm\xi,\label{pt}
\ee
where $\bm\theta_i$ was defined by (\ref{theta_i_def},  \ref{theta_i_unique}).
\end{definition}
\color{black}
Using the definition of  CPT's, one can easily show that
\begin{equation} \label{arbitraryshaped}
\ds \frac{1}{2} \sum^3_{i=1} \bH_0(\bz)_i \int_B{\bm D}^2_\bx
G(\bx,\bz)\bm\xi\times(\bm\theta_i+{\itbf e}_i\times\bm\xi)d\bm\xi
 =\ds  \sum^3_{l,l^\prime=1} {\bm D}^2_\bx G(\bx,\bz)_{ll^\prime} \mathbb{M}^{l,l^\prime} \bH_0(\bz). \end{equation}
\color{black}

Now we assume that $\bJ_0$ is a dipole source whose position is
denoted by $\bs$
\begin{equation} \label{jo}
 \bJ_0(\bx)= \rot \big( \bp\,
\delta(\bx, \bs) \big)  ,
\end{equation}
where $\delta(\cdot, \bs)$ is the Dirac  mass at $\bs$ and the
unit vector $\bp$ is the direction of the magnetic dipole.  In the absence of any target, the magnetic
field $\bH_0$ due to $\bJ_0(\bx)$ is given by
\begin{equation}
\label{ehon} \bH_0(\bx) = \rot \rot (\bp G(\bx, \bs))={\bm D}^2_{\bx}G(\bx,\bs)\bp, \quad \forall \; \bx\neq \bs.
\end{equation}

Assume for the sake of simplicity that $\mu_0=\mu_*$. Therefore,  by (\ref{arbitraryshaped}),
the asymptotic formula in Theorem \ref{thm-1} can be rewritten as follows.

\begin{cor} Assume that  $\bJ_0$ is a dipole source given by (\ref{jo}). Then, 
\begin{equation} \label{finalformula}
\bq \cdot (\bH_\alpha-\bH_0)(\bx) \simeq  \i k \alpha^5
\ds  \sum^3_{l,l^\prime=1} {\bm D}^2_\bx G(\bx,\bz)_{ll^\prime} \bq\cdot\mathbb{M}^{l,l^\prime} {\bm D}^2_{\bx}G(\bz,\bs)\bp
\end{equation}
for any unit vector $\bq$, where $\mathbb{M}^{l,l^\prime}$, defined by (\ref{pt}), are the CPTs associated with $B$.
\end{cor}
Note that, by following exactly the same arguments as in \cite{ACCGV}, we can prove that (\ref{finalformula}) is valid not only for $\nu$ of order one but also for $\nu$ much smaller than one.

Next, writing
$$
\mathbb{M}  = \Re e \, \mathbb{M}  + \i \Im m \,  \mathbb{M} ,
$$
we obtain
\begin{equation} \label{for_MSR_1}
\Re e \big( \bq \cdot (\bH_\alpha-\bH_0)(\bx) \big) \simeq - k
\alpha^5 \sum_{l,l^\prime=1}^3{\bm D}^2_\bx G(\bx,\bz)_{ll^\prime}\bq\cdot(\Im m \,\mathbb{M}^{l,l^\prime} ) \left({\bm D}^2_\bx G(\bz,\bs)\bp\right),
\end{equation}
and
\begin{equation} \label{for_MSR_2}
\Im m \big( \bq \cdot (\bH_\alpha-\bH_0)(\bx) \big) \simeq  k
\alpha^5 \sum_{l,l^\prime=1}^3{\bm D}^2_\bx G(\bx,\bz)_{ll^\prime}\bq\cdot(\Re e \,\mathbb{M}^{l,l^\prime} ) \left({\bm D}^2_\bx G(\bz,\bs)\bp\right).
\end{equation}

\color{black}
\begin{definition}
Let $\bs_m,m=1,2,\ldots,M$, be $M$ fixed points in $\RR^3$. These points
will be referred to as  sources.
Let $\bm{r}_n,n=1,2,\ldots,N$, be $N$ fixed points in $\RR^3$. These points
will be referred to as  receivers.
Fix two vectors $\bp$ and $\bq$ in $\RR^3$ and define the magnetic vector field
$\bH_0(\bx) = {\bm D}^2_{\bx}G(\bx,\bs_m)\bp$.
Define a perturbed field $\bE_\alpha$ as in (\ref{eform})
for the forcing term $\bJ_0(\bx)= \rot \big( \bp\,
\delta(\bx, \bs_m) \big) $, and set $\bH_\alpha= \rot\bE_\alpha / (i \omega \mu_\alpha)$.
Assume that all the receivers $\bm{r}_n $ and the sources 
$\bs_m$ are some positive distance away from the conductive object $\alpha B$ involved
in defining  $\bE_\alpha$.
We define the $MSR$ matrix $A$ to be the $N$ by $M$ matrix whose
 $n m$ -th entry is
\be \label{MSR_0}
A_{nm}=
\bq \cdot \Re e(\bH_\alpha-\bH_0)(\bm{r}_n).
\ee
\end{definition}
In the case where $\mu_\alpha$ is uniformly equal to $\mu_0$ and  $\alpha$ is small 
while $\nu$ defined in (\ref{deff5}) is $O(1)$, 
asymptotic formulas (\ref{for_MSR_1}) and (\ref{for_MSR_1}) lead to the estimate
for the $nm$ -th entry of the MSR matrix $A$
\be \label{MSR_def}
A_{nm}=k
\alpha^5 \sum_{l,l^\prime=1}^3{\bm D}^2_\bx G(\bm{r}_n,\bz)_{ll^\prime}\bq\cdot(\Re e \,\mathbb{M}^{l,l^\prime} ) \left({\bm D}^2_\bx G(\bz,\bs_m)\bp\right)+R_{nm},
\ee
where $R_{nm}$ is for lower order terms appearing due to the use of these asymptotic formulas.

\color{black}
\section{Statement of the \textsl{detection and identification problem} studied in this paper}
\label{new_section}
Let $I$ be a finite number and ${\cal{C}} = \{B^1, B^2, \ldots, B^I\}$ a collection of bounded domains
in $\RR^3$.
Let $B_\alpha$ be a domain obtained by dilation, rotation, and translation, of an element in ${\cal{C}} $:
$$
B_\alpha = \alpha R B^i + z,
$$
where $i$ is in $\{1,\ldots, I \}$, $\alpha >0 $, $R$ is a rotation, and $z$ is in $\RR^3$.
Assume that $B_\alpha$ has some (unknown) conductivity $\sigma >0$ and that using the eddy current
defined by (\ref{eform}) we can form the MSR matrix $A$ defined in
(\ref{MSR_0}). 
The \textsl{detection and identification problem} studied in this paper can now be simply formulated by
asking:\\
\begin{center}
Given $A$, find $i$.
\end{center}
Although this question may at first sight appear trivial,  a lot of issues arise in practice.
Is the solution unique? How will measurement noise affect the search for a solution? 
Since it is known 
that the computational cost of Newton like methods for inverse problems can be prohibitive, 
can we find a non iterative method which \textbf{avoids the trouble of solving forward problem}
(\ref{eform})?
The core contribution of our work is that thanks to a detailed analysis of how 
dilations, rotations, and translations affect the  MSR matrix $A$,
we are able to derive invariant quantities computed from $A$, which in turn
makes it possible to build a non iterative detection and identification algorithm.
In subsequent sections, we proceed to explain in details what these invariant quantities are, how this algorithm
was built, and how well it performs on simulated data.
\color{black}
\section{Localization algorithm}\label{Localization algorithm}
Assume that 
measurements used in building the MSR matrix $A$
are  tinted by  noise.
In this paper we utilize
  Hadamard's sampling technique as proposed in \cite{ACCGV}:
this  is a data
acquisition scheme deigned to reduce  noise. It allows us to acquire simultaneously
all the elements of the MSR matrix while  reducing the effects of noise. The
main advantage to using  Hadamard's technique  is that it divides the variance of  measurement noise
by the number of sources \cite{AGJKLSW}.

Doing so, we can rewrite the MSR matrix in the following form
\be \label{rewrite A}
A=  \color{black}U\mathcal{M}_q V_p+R+ \frac{\sigma_{\mathrm{noise}}}{\sqrt{M}} W,
\label{withR}
\ee
where $R$ is a higher-order error term due to  using the asymptotic formula from
Theorem \ref{thm-1},
$W$ is a $N\times M$ matrix with independent and identical Gaussian entries with zero mean and unit variance, and
$\sigma_{\mathrm{noise}}$ is a small positive constant. The matrix  $U$ is a N-by-9 matrix of the form
\ben
U=\begin{pmatrix}
{\bm D}^2_\bx G(\br_1,\bz)_{11} & {\bm D}^2_\bx G(\br_1,\bz)_{12}& \ldots & {\bm D}^2_\bx G(\br_1,\bz)_{33}\\
\vdots& \vdots& \vdots& \vdots\\
{\bm D}^2_\bx G(\br_N,\bz)_{11}& {\bm D}^2_\bx G(\br_N,\bz)_{12}& \ldots &{\bm D}^2_\bx G(\br_N,\bz)_{33}
\end{pmatrix},
\een
$\mathcal{M}_q$ is a 9-by-3 matrix of the form
\be \label{Mq_def}
\mathcal{M}_q =k\alpha^5\Re e \begin{pmatrix}
\bq^T\mathbb{M}^{1,1}\\
\bq^T\mathbb{M}^{1,2}\\
\vdots\\
\bq^T\mathbb{M}^{3,3}
\end{pmatrix},
\ee
and $V_p$ is a 3-by-M matrix of the form
\ben
V_p=\begin{pmatrix}
 {\bm D}^2_\bx G(\bz,\bs_1)\bp &  {\bm D}^2_\bx G(\bz,\bs_2)\bp & \ldots &  {\bm D}^2_\bx G(\bz,\bs_M)\bp
\end{pmatrix}.
\een
Define the linear operator $L:\mathbb{R}^{9\times 3}\rightarrow \mathbb{R}^{N\times M}$ by
\be \label{L_def}
L(\mathcal{M}_q)=U\mathcal{M}_qV_p .
\ee
Dropping the lower-order term $R$ in (\ref{rewrite A}), the MSR matrix
can be approximated as follows
 $$A\approx  L(\mathcal{M}_q)+ \frac{\sigma_{\mathrm{noise}}}{\sqrt{M}} W.$$
If the target $B$ is a sphere, the operator can be simplified as $L(\mathcal{M}_q)=\mathcal{M}V_q^\prime V_p$, where $\mathcal{M}$ is a real scalar and $V_q$ is defined as $V_p$ with $\bq$ instead of $\bp$ (see \cite{ACCGV}).   We used the MUSIC algorithm to localize the spherical target.  In the present paper, for arbitrary shaped targets, let ${\bm P}$ be the orthogonal projection onto the right null space  of  $L(\mathcal{M}_q)$. We define the imaging functional as
\be
\mathcal{I}_{MU}(\bz^S)=\big[\frac{1}{\sum^3_{i=1}\|{\bm P}({\bm D}^2_\bx G(\bz^S,\bs_1)\bp\cdot\bm{e}_i, {\bm D}^2_\bx G(\bz^S,\bs_2)\bp\cdot \bm{e}_i, \ldots, {\bm D}^2_\bx G(\bz^S,\bs_M)\bp\cdot\bm{e}_i)\|^2}\big]^{1/2} \label{music}
\ee
for $\bz^S$ in the search domain.
Following \cite{AILP},  we obtain the following result.
\begin{prop}
Suppose that $U \mathcal{M}_q$ has full rank. Then $L(\mathcal{M}_q)$ has three non zero singular values. Furthermore,
$\mathcal{I}_{MU}(\bz^S)$ attains its maximum approximately at $\bz^S=\bz$.
\end{prop}

As it will be shown in section 6, the MUSIC algorithm still works for arbitrary shaped targets. In section 6, we also numerically investigate the resolution of the MUSIC imaging algorithm in the presence of measurement noise.
\section{Properties of the CPTs $\mathbb{M}^{l,l^\prime}$}\label{Properties of the CPTs}

We call a {dictionary} a
collection of standard shapes, which are centered at the origin
and with characteristic sizes of order 1. Given the CPTs of an
unknown shape $D$, and assuming that $D$ is obtained from a
certain element $B$ in the dictionary by applying some unknown rotation
$\theta$, scaling $s$ and translation $\bz$, our objective is to recognize $B$ from the dictionary using induction data at a single or multiple frequencies.
For doing so,
one may proceed by first reconstructing the shape $D$ using its
CPTs through some optimization procedures, and then match the reconstructed shape with
the dictionary. However, such a method may be time-consuming and the
recognition efficiency depends on the shape reconstruction
algorithm.

We propose a shape identification algorithm
using the CPTs. The algorithm operates directly in
the data domain which consists of CPTs and avoid the need for
reconstructing the shape $D$. The heart of our approach is some
invariance relations between the CPTs of $D$ and $B$.

We first establish the following lemma.

\begin{lem}\label{lem-31}
Let ${\cal{O}}$ be an orthogonal $3\times 3$ matrix. 
\begin{itemize}
\item[{\rm (i)}] If $\bu, \bv$ are two vectors in $\RR^3$
 then
\be
 ({\cal{O}} \bu \times {\cal{O}}\bv) = (\det {\cal{O}}) {\cal{O}} (\bu \times \bv), \\
  \bu \times ({\cal{O}}\bv) = (\det {\cal{O}}) {\cal{O}} ( ({\cal{O}}^T \bu) \times \bv).
\ee
\item[{\rm (ii)}] If $\bF$ is a $\mathcal{C}^1$-vector field in $\RR^3$ then
\be
\rot ({\cal{O}}^T \bF({\cal{O}} \bx)) = (\det {\cal{O}} ) {\cal{O}}^T (\rot \bF ) ({\cal{O}} \bx),
\label{rotinv}\\
\rot \rot ({\cal{O}}^T \bF({\cal{O}} \bx)) = {\cal{O}}^T (\rot \rot \bF ) ({\cal{O}} \bx),
\label{rotinv2}\\
\ddiv ({\cal{O}}^T \bF({\cal{O}} \bx)) = (\ddiv \bF) ({\cal{O}} \bx) .\label{divinv}
\ee
\item[{\rm (iii)}] If $\bn$ is the outward normal vector on a $\mathcal{C}^1$-surface which is invariant
under ${\cal{O}}$ then
\be
 \bn({\cal{O}}\bx) = {\cal{O}} \bn(\bx).
\ee
\end{itemize}
\end{lem}
\debproof
 (i) is due to the fact that ${\cal{O}}$ maps orthonormal basis to
orthonormal basis. \\
Formula (\ref{rotinv}) is most easily shown by Fourier transform.
Without loss of generality we may assume that $\bF$ has compact support. We first note that  if $\bG$ is any compactly supported  $\mathcal{C}^1$-vector field in $\RR^3$,
then
\ben
\widehat{\rot \bG} (\bm \zeta) =  \int  e^{\i \bx \cdot \bm \zeta}
\rot \bG(\bx) d \bx
= - \i  \int  e^{\i \bx \cdot\bm \zeta}
\bm\zeta \times \bG(\bx) d \bx,
\een
and
\ben
\widehat{\bG} ({\cal{O}} \bm\zeta)  =  \int  e^{\i \bx \cdot {\cal{O}} \bm\zeta}
 \bG(\bx) d \bx
=  \int  e^{\i {\cal{O}}\bx \cdot {\cal{O}}\bm\zeta}
 \bG({\cal{O}} \bx) d \bx
= \int  e^{\i \bx \cdot  \bm\zeta}
 \bG({\cal{O}} \bx) d \bx
= \widehat{\bG({\cal{O}} \bx)} ( \bm\zeta).
\een
Using these two formulas and the notation ${\cal{F}}$ for Fourier transforms we write
\ben
{\cal{F}} (  \nabla \times ({\mathcal{O}}^T\bF({\mathcal{O}} \bx))) (\bm\zeta) =
 \int  e^{\i \bx \cdot \bm\zeta}
\rot ( {\cal{O}}^T\bF({\cal{O}} \bx) ) d \bx
= - \i  \int  e^{\i \bx \cdot\bm \zeta}
\bm\zeta \times ( {\cal{O}}^T\bF({\cal{O}} \bx) ) d \bx \\
=-\i (\det {\cal{O}}) {\cal{O}}^T  \int  e^{\i {\cal{O}}\bx \cdot {\cal{O}} \bm\zeta} ({\cal{O}} \bm\zeta) \times
\bF({\cal{O}} \bx)  d \bx
= -\i (\det {\cal{O}}) {\cal{O}}^T  \int  e^{\i \bx \cdot {\cal{O}}\bm\zeta} ({\cal{O}} \bm\zeta) \times
\bF(\bx)  d \bx \\
= (\det {\cal{O}}) {\cal{O}}^T \widehat{\rot \bF} ({\cal{O}} \bm\zeta),
\een
which yields formula (\ref{rotinv}). \\
Formula (\ref{rotinv2}) follows easily from  (\ref{rotinv}) and formula
 (\ref{divinv}) is proved likewise. \\
To prove (iii) we can assume that the surface is given by the equation
$f(\bx) = 0$, where $f$ satisfies $f({\cal{O}} \bx) = f(\bx)$.
It follows that ${\cal{O}}^T (\nabla f)({\cal{O}} \bx) = \nabla f(\bx )$ so
$ (\nabla f)({\cal{O}} \bx) ={\cal{O}} (\nabla f(\bx ))$ and $ \|(\nabla f)({\cal{O}} \bx)\|
 =\| \nabla f(\bx )\|$ and  formula in (iii) holds.
\finproof

Let $B_\bz=\bz+B$ be a shift of $B$. Denote $\mathbb{M}^{l,l^\prime}_i[B_\bz]$ be the $i$-th column of the conductive polarization tensor. The following result holds.
\begin{prop}[translation formula] \label{prop-shift}
$\mathbb{M}^{l,l^\prime}_i[B_\bz]=\mathbb{M}^{l,l^\prime}_i[B]$.
\end{prop}

\debproof Let $F_{\bz}$ be the solution to the problem
\ben
\nabla_\xi \times\mu^{-1}\nabla_\xi \times F_{\bz}-\i\omega\sigma\alpha^2F_{\bz}=\i\omega\sigma\alpha^2\bm{e}_i\times\bm\xi \mbox{ in } B_\bz\cup {B_\bz}^c,\\
\nabla_\xi \cdot F_{\bz}=0 \mbox{ in } B^c_\bz,\\
\mbox{[}\bn\times F_{\bz}\mbox{]}=0 \mbox{ on } \partial B_\bz,\\
\mbox{[}\mu^{-1}\nabla_\xi\times F_{\bz}\times \bn] =-2[\mu^{-1}]\bm{e}_i\times\bn \mbox{ on } \partial B_\bz,\\
\color{black}\int_{\partial B_{\bz}} (F_{\bz}\cdot \bn)^+ =0, \color{black}\\
F_{\bz}=O(|\bm\xi|^{-1}) \mbox{ as } |\bm\xi|\rightarrow \infty .
\een
Define $F_0$ to be equal to $F_{\bz}$ for the choice $\bz=0$.
It can be easily seen  that
$$F_{\bz}=F_0+G_\bz,$$
where $G_\bz$ solves
\ben
\nabla_\xi \times\mu^{-1}\nabla_\xi\times G_\bz-\i\omega\sigma\alpha^2 G_\bz = \i\omega\sigma\alpha^2\bm{e}_i\times\bz \mbox{ in } B\cup {B}^c,\\
\nabla_\xi\cdot G_\bz=0 \mbox{ in } B^c,\\
\mbox{[}\bn\times G_\bz]=0, \mbox{[}\mu^{-1}\nabla_\xi\times G_\bz\times\bn]=0 \mbox{ on }\partial B,\\
\int_{\partial B} (G_{\bz}\cdot \bn)^+ =0, \\
G_\bz=O(|\bm\xi|^{-1}) \mbox{ as } |\bm\xi|\rightarrow\infty.
\een
Let $\nabla u = -\bm{e}_i\times \bz$. Then, due to the fact that $\bm{e}_i\times\bz$ is a constant vector, $u$ is a linear function.
Let  $\tilde u$  be defined by
\ben
\left\{\begin{array}{c}
\Delta\tilde{u}=0 \mbox{ in } B^c,\\
\nm
\tilde{u} =u \mbox{ on } \partial B,\\
\nm
\tilde{u}=O(|\bm\xi|^{-1})\mbox{ as } |\bm\xi|\rightarrow\infty.
\end{array}\right.
\een
We have thus determined $G_\bz$. It can be expressed as 
$$G_\bz=\left\{\begin{array}{c} -\bm{e}_i\times\bz\mbox{ in }
 B,\\ \nm \nabla\tilde{u}\mbox{ in } {B^c}. \end{array}\right.$$
Note that  $\nabla\times\nabla\times G_\bz=0$ in $B$. Therefore, it follows that the $i$-th column of $\mathbb{M}^{l,l^\prime}$ is given by
\ben
\mathbb{M}^{l,l^\prime}_i[B_\bz]&=&\frac 12 \bm{e}_l\times\int_{B_\bz}\xi_{l^\prime}(F_{\bz}+\bm{e}_i\times\bm\xi)d\bm\xi\\
&=&\frac 12 \frac{1}{\i\omega\sigma\alpha^2\mu}\bm{e}_l\times\int_{B_\bz}\xi_{l^\prime}(\nabla\times\nabla\times F_{\bz})d\bm\xi\\
&=&\frac 12\frac 1{\i\omega\sigma\alpha^2\mu}\bm e_l\times\int_{B}(\bm z_{l'}+\hat{\bm\xi}_{l'})(\nabla_{\hat{\bm\xi}}\times\nabla_{\hat{\bm\xi}}\times F_{\bm z}(\bm z+\hat{\bm\xi}))d\hat{\bm\xi}\\
&=&\mathbb{M}^{l,l'}_i[B]+\frac 12\frac 1{\i\omega\sigma\alpha^2\mu}\bm z_{l'}\int_B\bm e_l\times(\nabla\times\nabla\times\bm\theta_i)d\bm\xi\\
&=&\mathbb{M}^{l,l^\prime}_i[B].
\een
In the last equality, we have used the fact that $\int_B\bm e_l\times(\nabla\times\nabla\times\bm\theta_i)d\bm\xi = 0$ which
is proved in \cite{ACCGV}.
\finproof

Let $s>0$ be a scaling factor. Let $sB$ be the scaled domain and let $\mathbb{M}^{l,l^\prime}[\omega\sigma,sB]$ be the conductive polarization tensor associated with the scaled domain $sB$.
\begin{prop}[scaling formula] \label{prop-scaling} We have the following scaling relation:
$$\mathbb{M}^{l,l^\prime}[\omega\sigma,sB]=s^5\mathbb{M}^{l,l^\prime}[\omega\sigma s^2,B].$$
\end{prop}
\debproof
Let $F_{\omega\sigma,sB}(\bm\xi)$ be defined by the interface problem
\ben
\nabla\times\mu^{-1}\nabla\times F_{\omega\sigma,sB}(\bm\xi)-\i\omega\sigma\alpha^2
F_{\omega\sigma,sB}(\bm\xi)&=&\i\omega\sigma\alpha^2\bm{e}_i\times\bm\xi \mbox{ in } (sB)\cup {(sB)}^c,\\
\nabla\cdot F_{\omega\sigma,sB}&=&0\mbox{ in } (sB)^c,\\
\mbox{[}\bn\times F_{\omega\sigma,sB}]=0,\mbox{[}\mu^{-1}\nabla\times F_{\omega\sigma,sB}\times\bn]&=&-2[\mu^{-1}]\bm{e}_i\times\bn\mbox{ on }\partial (sB),\\
\int_{\partial s B} (F_{\omega\sigma,sB}\cdot \bn)^+ &=&0, \\
F_{\omega\sigma,sB}(\bm\xi)&=&O(|\bm\xi|) \mbox{ as } |\bm\xi|\rightarrow\infty,
\een
where all the gradients are taken in the $\xi$ variable.
Setting $\bm\xi = s\bm\xi'$, it follows that
\ben
\nabla\times\mu^{-1}\nabla\times ( F_{\omega\sigma,sB}(s\bm\xi'))-\i\omega\sigma\alpha^2s^2F_{\omega\sigma,sB}(s\bm\xi')&=&\i\omega\sigma\alpha^2s^2\bm{e}_i\times s\bm\xi' \mbox{ in } B\cup {B}^c,\\
\nabla\cdot F_{\omega\sigma,sB}(s\bm\xi')&=&0\mbox{ in } B^c,\\
\mbox{[}\bn\times F_{\omega\sigma,sB}]=0,\mbox{[}\mu^{-1}\nabla\times (F_{\omega\sigma,sB}(s\bm\xi'))\times\bn]&=&-2s[\mu^{-1}]\bm{e}_i\times\bn\mbox{ on }\partial B,\\
\int_{\partial B} (F_{\omega\sigma,sB} (s\bm\xi')\cdot \bn)^+ d \bm\xi' &=&0, \\
F_{\omega\sigma,sB}(s\bm\xi')&=&O(|\bm\xi'|) \mbox{ as } |\bm\xi'|\rightarrow\infty.
\een
These equations indicate that
$$\frac 1s F_{\omega\sigma,sB}(s\bm\xi)=F_{\omega\sigma s^2,B}(\bm\xi).$$
Hence,
\ben
\mathbb{M}^{l,l^\prime}_i[\omega\sigma,sB]&=&\frac 12\bm{e}_l\times\int_{sB}\xi_{l^\prime}(F_{\omega\sigma,sB}+\bm{e}_i\times\bm\xi)d\bm\xi\\
&=&\frac 12s^3\bm{e}_l\times\int_B(s\xi')_{l^\prime}(F_{\omega\sigma,sB}(s\xi')+\bm{e}_i\times(s\bm\xi'))d\bm\xi'\\
&=&\frac 12 s^5\bm{e}_l \times \int_B\xi_{l^\prime}(F_{\omega\sigma s^2,B}(\xi)+\bm{e}_i\times\xi)d\bm\xi\\
&=&s^5\mathbb{M}^{l,l^\prime}[\omega\sigma s^2,B] ,
\een
which completes the proof.
\finproof

Let ${\cal{O}}$ be a rotation of $\mathbb{R}^3$ whose axis passes through the origin. We also denote by ${\cal{O}}$ its matrix in the 
natural basis of  $\mathbb{R}^3$. 
Let $\mathbb{M}^{l,l^\prime}$ be the conductive polarization tensor associated with the domain $B$. 
It proves convenient to reshape 
 the 9 CPT matrices for the domain $B$  as follows:
$$\mathbb{M}[\omega\sigma, B]=\begin{pmatrix}\mathbb{M}^{1,1}&\mathbb{M}^{1,2}&\mathbb{M}^{1,3}\\
\mathbb{M}^{2,1}&\mathbb{M}^{2,2}&\mathbb{M}^{2,3}\\
\mathbb{M}^{3,1}&\mathbb{M}^{3,2}&\mathbb{M}^{3,3}\\
\end{pmatrix},$$  and to denote by $\mathbb{M}[\omega\sigma, {\cal{O}}(B)]$ its  counterpart relative
to the rotated domain ${\cal{O}}(B)$.
 We obtained the following result.
\begin{prop}[rotation formula] \label{prop-rotation} 
The following identity holds
$$\mathbb{M}[\omega\sigma, {\cal{O}}(B)]={\cal{O}}_2{\cal{O}}_1\mathbb{M}[\omega\sigma,B]{\cal{O}}_1^T{\cal{O}}_2^T,$$
where ${\cal{O}}_1$ is the $9 \times 9$  matrix defined by the blocks $diag({\cal{O}},{\cal{O}},{\cal{O}})$,
${\cal{O}}_2$ is the $9 \times 9$   matrix
defined by the blocks
$\begin{pmatrix}
{\cal{O}}_{11}I_3&{\cal{O}}_{12}I_3&{\cal{O}}_{13}I_3\\
{\cal{O}}_{21}I_3&{\cal{O}}_{22}I_3&{\cal{O}}_{23}I_3\\
{\cal{O}}_{31}I_3&{\cal{O}}_{32}I_3&{\cal{O}}_{33}I_3\\
\end{pmatrix}$,
and ${\cal{O}}_{ij}$ is the ij-th entry of ${\cal{O}}$.
\end{prop}
\debproof
Denote by $F_{{\cal{O}}(B),\bm{e}_i}(\bm\xi)$  the solution to the interface problem
\ben
\nabla\times\mu^{-1}\nabla\times F_{{\cal{O}}(B),\bm{e}_i}-\i\omega\sigma\alpha^2 F_{{\cal{O}}(B),\bm{e}_i}&=&\i\omega\sigma\alpha^2\bm{e}_i\times\bm\xi \mbox{ in } {\cal{O}}(B)
\cup {{\cal{O}}(B)}^c,\\
\nabla\cdot F_{{\cal{O}}(B),\bm{e}_i}&=&0 \mbox{ in } {\cal{O}}(B)^c,\\
\mbox{[}\bn \times F_{{\cal{O}}(B),\bm{e}_i}\mbox{]}&=&0 \mbox{ on } \partial {\cal{O}}(B),\\
\mbox{[}\mu^{-1}\nabla\times F_{{\cal{O}}(B),\bm{e}_i}\times\bn\mbox{]}&=&-2\mbox{[}\mu^{-1}\mbox{]}\bm{e}_i\times\bn \mbox{ on } \partial {\cal{O}}(B),\\
 \int_{\partial {\cal{O}}(B)}        (F_{{\cal{O}}(B),\bm{e}_i} \cdot n)^+ &=&0              ,\\
F_{{\cal{O}}(B),\bm{e}_i}(\bm\xi)&=&O(|\bm\xi|^{-1}) \mbox{ as } |\bm\xi| \rightarrow \infty.
\een
Next, we apply identities from Lemma \ref{lem-31} to obtain
\ben
{\cal{O}}\nabla\times\mu^{-1}\nabla\times ( {\cal{O}}^T F_{{\cal{O}}(B),\bm{e}_i}({\cal{O}}\bm\xi) )-\i\omega\sigma\alpha^2 F_{{\cal{O}}(B),\bm{e}_i}({\cal{O}}\bm\xi) = \i\omega\sigma\alpha^2\bm{e}_i\times {\cal{O}}\bm\xi
\mbox{ in } B
B^c,
\een
\ben
\nabla\cdot({\cal{O}}^T F_{{\cal{O}}(B),\bm{e}_i}({\cal{O}}\xi))=0 \mbox{ in } B^c,
\een
\ben
\mbox{[} (\bn(\bm\xi)\times {\cal{O}}^T F_{{\cal{O}}(B),\bm{e}_i}({\cal{O}}\bm\xi) \mbox{]} =0 
\mbox{ on } \partial B,
\een
 and
\ben
\mbox{[}\mu^{-1}{\cal{O}} \nabla\times ({\cal{O}}^T  F_{{\cal{O}}(B),\bm{e}_i}({\cal{O}}\xi))\times\bn({\cal{O}}\bm\xi)\mbox{]}
= {\cal{O}} \mbox{[} \mu^{-1} \nabla\times ({\cal{O}}^T  F_{{\cal{O}}(B),\bm{e}_i}({\cal{O}}\xi)) \times\bn(\bm\xi) \mbox{]} \\
= {\cal{O}}(-2\mbox{[}\mu^{-1}\mbox{]} {\cal{O}}^T\bm{e}_i\times\bn(\xi)) \mbox{ on } \partial B.
\een
Thus we get the relation:
$$F_{{\cal{O}}(B),\bm{e}_i}({\cal{O}}\bm\xi)={\cal{O}}F_{B,{\cal{O}}^T\bm{e}_i}(\bm\xi), \; \forall\bm\xi\in \mathbb{R}^3.$$

Now, using the definition of the conductive polarization tensor, we obtain that
\ben
\mathbb{M}^{l,l^\prime}_i[{{\cal{O}}}(B)]&=&\frac 12 \bm{e}_l\times\int_{{{\cal{O}}}(B)}\xi_{l^\prime}(F_{{{\cal{O}}}(B),\bm{e}_i}(\bm\xi)+\bm{e}_i\times\bm\xi)d\bm\xi\\
&=&\frac 12\bm{e}_l\times\int_B({\cal{O}}\bm\xi)_{l^\prime}({\cal{O}}F_{B,{\cal{O}}^T\bm{e}_i}(\bm\xi)+\bm{e}_i\times 
{\cal{O}}\bm\xi)d\bm\xi\\
&=&\frac 12 \bm{e}_l\times {\cal{O}}\int_B({\cal{O}}\bm\xi)_{l^\prime}(F_{B,{\cal{O}}^T\bm{e}_i}(\bm\xi)+{\cal{O}}^T\bm{e}_i\times\bm\xi)d\bm\xi\\
&=&\frac 12 {\cal{O}}\int_B({\cal{O}}\bm\xi)_{l^\prime}({\cal{O}}^T\bm{e}_l)\times (F_{B,{\cal{O}}^T\bm{e}_i}(\bm\xi)+{\cal{O}}^T\bm{e}_i\times\bm\xi)d\bm\xi\\
&=&\frac 12 {\cal{O}}\sum^3_{m=1}{\cal{O}}_{lm}\int_B({\cal{O}}\bm\xi)_{l^\prime}\bm{e}_m\times(F_{B,{\cal{O}}^T\bm{e}_i}(\bm\xi)+{\cal{O}}^T\bm{e}_i\times\bm\xi)d\bm\xi\\
&=&\frac 12 {\cal{O}}\sum^3_{m=1}{\cal{O}}_{lm}\int_B\sum^3_{n=1}{\cal{O}}_{l^\prime n}\xi_n\bm{e}_m\times(F_{B,{\cal{O}}^T\bm{e}_i}(\bm\xi)+{\cal{O}}^T\bm{e}_i\times\bm\xi)d\bm\xi\\
&=&{\cal{O}}\sum^3_{m=1}{\cal{O}}_{lm}\sum^3_{n=1}{\cal{O}}_{l^\prime n}\mathbb{M}^{m,n}_i[B,{\cal{O}}^Te_i].
\een
From ${\cal{O}}^Te_i\times\bm\xi=\sum^3_{p=1}{\cal{O}}_{ip}\bm{e}_p\times\bm\xi$, we have
$$F_{B,{\cal{O}}^Te_i}=\sum^3_{p=1}{\cal{O}}_{ip}F_{B,\bm{e}_p}$$ and
\ben
\frac 12\bm{e}_m\times\int_B\xi_n(F_{B,{\cal{O}}^T\bm{e}_i}(\bm\xi)+{\cal{O}}^T\bm{e}_i\times\bm\xi)d\bm\xi\\
=\frac 12\bm{e}_m\times\int_B\xi_n\sum^3_{p=1}{\cal{O}}_{ip}(F_{B,\bm{e}_p}(\bm\xi)+\bm{e}_p\times\bm\xi)d\bm\xi
=\sum^3_{p=1}{\cal{O}}_{ip}\mathbb{M}^{m,n}_p.
\een
Finally, we arrive at
$$\mathbb{M}^{l,l^\prime}_i[{{\cal{O}}}(B)]={\cal{O}}\sum^3_{m=1}{\cal{O}}_{lm}\sum^3_{n=1}{\cal{O}}_{l^\prime n}\sum^3_{p=1}{\cal{O}}_{ip}\mathbb{M}^{m,n}_p[B].$$
At this stage we observe that
$$
\sum^3_{p=1}{\cal{O}}_{ip}\mathbb{M}^{m,n}_p[B] = (\mathbb{M}{\cal{O}}_1^T)^{m,n}_i ,
$$
and further simplifications lead to
the desired result. \finproof
Recall that ${\cal{O}}$ is a unitary matrix. In view of the special structure of ${\cal{O}}_1,{\cal{O}}_2$, we have the following results.
\begin{prop}
The matrices ${\cal{O}}_1$ and ${\cal{O}}_2$ are orthogonal matrices. Moreover, $\mathbb{M}[\omega\sigma,{{\cal{O}}}(B)]$ and $\mathbb{M}[\omega\sigma,B]$ have the same singular values.
\end{prop}
\begin{rem}
Proposition \ref{prop-rotation} expresses the fact that the singular values of $\mathbb{M}$ are invariant under rotations,
and Proposition \ref{prop-shift} that CPT's are invariant under translations.
Consequently, the singular values of $\mathbb{M}$ are invariant under translations and rotations.
As to Proposition \ref{prop-scaling}, it indicates how  CPTs depend on frequency.
\end{rem}
\begin{rem}\label{low-approx}
Since $F_{\omega\sigma s^2, B}$ solves the problem
\ben
\nabla\times\mu^{-1}\nabla\times F_{\omega\sigma s^2, B}(\bm\xi)-\i\omega\sigma\alpha^2s^2F_{\omega\sigma s^2,B}(\bm\xi)&=&\i\omega\sigma\alpha^2s^2\bm{e}_i\times\bm\xi \mbox{ in } B\cup {B}^c,\\
\nabla\cdot F_{\omega\sigma s^2, B}&=&0 \mbox{ in } B^c
\een
with boundary and interface conditions independent of $s$, 
it is clear that this vector field is continuous in $s$.
It follows that $\mathbb{M}[\omega\sigma s^2,B]$
and thanks to Proposition \ref{prop-scaling},
  $\mathbb{M}[\omega\sigma,sB]$, are also 
continuous in $s$. 

\end{rem}

\section{CPTs recovery and dictionary matching} \label{CPTs recovery}
\subsection{CPTs recovery}
Recall the definition of $\mathcal{M}_q$ given in equation (\ref{Mq_def}).
An approximation to the projection
of $\mathcal{M}_q$ on the orthogonal of the nullspace of linear operator $L$ defined in
(\ref{L_def})
can be formed by solving the following least squares minimization problem
\be
\mathcal{M}_q=\mathrm{arg} \min_{\mathcal{M}_q \bot \mathrm{ker}(L)}\|A-L(\mathcal{M}_q)\|^2_F ,\label{lsq-msr0}
\ee
where the MSR matrix $A$ was defined in  (\ref{MSR_0}) (a useful approximate
was given in (\ref{MSR_def})) and
$\|\cdot\|_F$ denotes the Frobenius norm of matrices.
It is clear that $L(\mathcal{M}_q)$ only depends on $\mathcal{M}_q$, the product of vector $\bq$ and CPT matrices and a scaling factor $k\alpha^5$.  So, for a given $\bq$, we can not recover all the entries of the CPT matrices.  If we let $\bm{q}=\bm{e}_1,\bm{e}_2, \bm{e}_3$, respectively, and solve the least squares problem three times,  we can recover all the entries of the CPT up to a scaling factor.   Furthermore, by \eqref{pt}, we know that  the entries in $l$-th row of  $\mathbb{M}^{l,l^\prime}$ are zeros,  we should find a solution  of \eqref{lsq-msr0} such that the $l$-th row of $\mathbb{M}^{l,l^\prime}$ is zero vector.

More precisely, for $\bq=\bm{e}_l$,  let $$M(\bm{e}_l)=\bigg\{M\in \mathbb{R}^{9\times 3}: \,  M_{ij}=0, i=3(l-1)+1, 3(l-1)+2,3(l-1)+3, j=1,2,3 \bigg\}.$$  The solution to the least squares problem
\be
\mathcal{M}_{\bm{e}_l}= \mathrm{arg} \min_{\mathcal{M}_q \in \mathrm{ker}(L)^\bot\cap M(\bm{e}_l)}\|A-L(\mathcal{M}_q)\|^2_F\label{lsq-msr}
\ee
for $l=1,2,3,$ will give the desired reconstruction of projections of CPT matrices on $\mathrm{ker}(L)^\bot$. 

\subsection{Dictionary matching}
 The CPT matrix  depends non linearly  on the scaling factor $s$ as shown in Proposition \ref{prop-scaling}.
Moreover,  recovery of CPT's can only be done  up to the scaling factor $k\alpha^5$. Consequently,
 we can not build the dictionary directly from the singular values of matrix $\mathbb{M}[\omega\sigma,B]$.
 Instead, we use the normalized singular values at multiple frequencies as the elements of the dictionary.
 More specifically,  we build the dictionary from the singular values of $\mathbb{M}$ for  multiple frequencies
 $\omega_n,n=1,2, \ldots, F$. In other words,  if we denote by $S^i_n$ the singular values of $\mathbb{M}$ for shape
 $i$ at frequency $\omega_n$ for $n=1,2,\ldots,F$, the corresponding element for this target in the dictionary is
$$D_i=[S^i_1,S^i_2,\ldots,S^i_F]/\max\{S^i_n,n=1,2,\ldots,F\},$$
and the dictionary is
$$\mathcal{D}=\{D_1, D_2, \ldots, D_I\},$$
where $I$ is the number of shapes in the dictionary.
This motivates us to implement the following dictionary matching algorithm.
\begin{alg}\label{alg-1}
Given the 
MSR matrices for $\bq=\bm{e}_1,\bm{e}_2,\bm{e}_3$ at frequency $\omega_n,n=1,2,\ldots,F$.
\begin{enumerate}
\item [Step 1.] At each frequency $\omega_n$, recover the CPT matrix $\mathbb{M}^{l,l^\prime}$ by successively 
setting $\bq=\bm{e}_1,\bm{e}_2,\bm{e}_3$  in \eqref{lsq-msr}, and forming the corresponding matrix $\mathbb{M}$.
\item [Step 2.] Apply the Singular Value Decomposition to $\mathbb{M}$ at each frequency $\omega_n$ and form  the vector $\hat{D}=[\hat{S}_1,\hat{S}_2,\ldots,\hat{S}_F]/\max_n\{\hat{S}_n\} $.
\item [Step 3.] Find the closest match to $\hat{D}$
within the dictionary $\mathcal{D}$ of precomputed elements $D$ by solving the minimization problem
              $\min_{D\in \mathcal{D}}\{\|D-\hat{D}\|_2\}.$ This
              will determine the approximate  shape of the target.
\end{enumerate}
\end{alg}
\section{Numerical examples}\label{Numerical examples}

\subsection{Testing  the MUSIC- type localization algorithm}
\begin{figure}
	\begin{center}
\includegraphics[scale=.5]{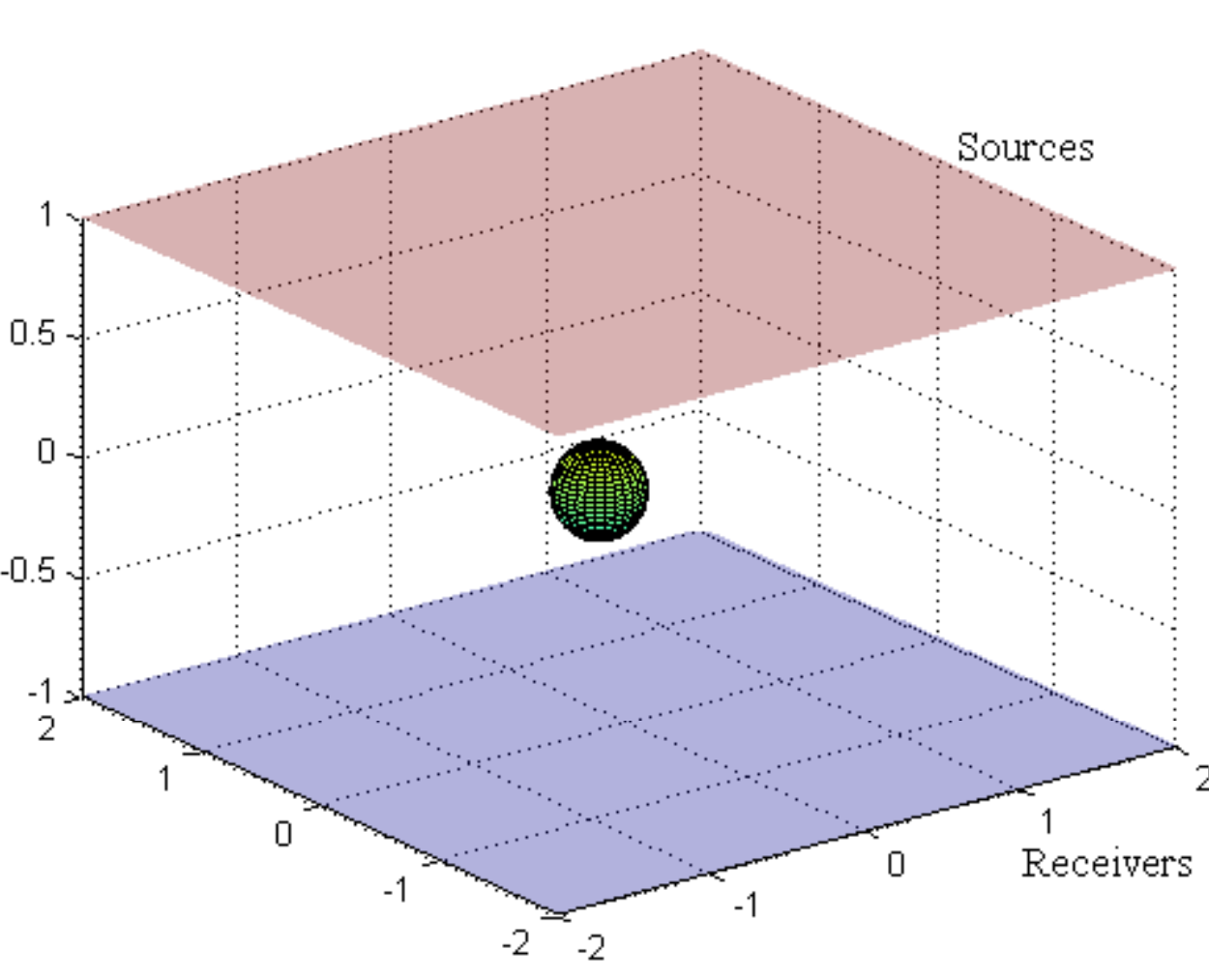}
\end{center}
\caption{\color{black} Sketch of the geometry for numerical simulations of 
shape detection and classification. 256 sources are   placed on a uniform grid on the top square $[-2,2]^2\times\{1\}$
and 256 receivers are   placed on a uniform grid on the bottom square $[-2,2]^2\times\{-1\}$.
An unknown conductive object that has to be detected and classified lies between these two squares.
}\label{geometry}
\end{figure}

We first 
illustrate how well our MUSIC- type algorithm 
performs the task of locating targets that are not necessarily spherical.
Pick   an ellipsoid shaped target with equation $x^2+y^2+\frac{z^2}{4}\leq \alpha^2$. 
The number of sources $M $ and the number of receivers $N$ are both chosen to be equal to
256 and are placed as indicated in Figure \ref{geometry}. 
We set the values 
$\sigma=5.97e7$S/m, $\mu_*=\mu_0= 1.2566e-06$ H/m, and $\alpha=0.01$ m, $\omega=133.5$, so that 
$k\alpha^2=1$ and the asymptotic formula form Theorem \ref{thm-1}  is valid.
Relevant CPTs are computed by a finite 
element code based on PHG \cite{phg}, 
and we then form the product $U\mathcal{M}_q V_p$.
To simulate the matrix $A$ (recall formula (\ref{withR})),
we generate an $N$ by $M$ matrix $W$ with entries from a normal distribution with mean 0 and variance 1
using the matlab function 'randn' and we compute the sum
$U\mathcal{M}_q V_p+  \frac{\sigma_{\mathrm{noise}}}{\sqrt{M}} W$ for different  values of $\sigma_{\mathrm{noise}} $.

\begin{figure}
\includegraphics[width=0.6\textwidth]{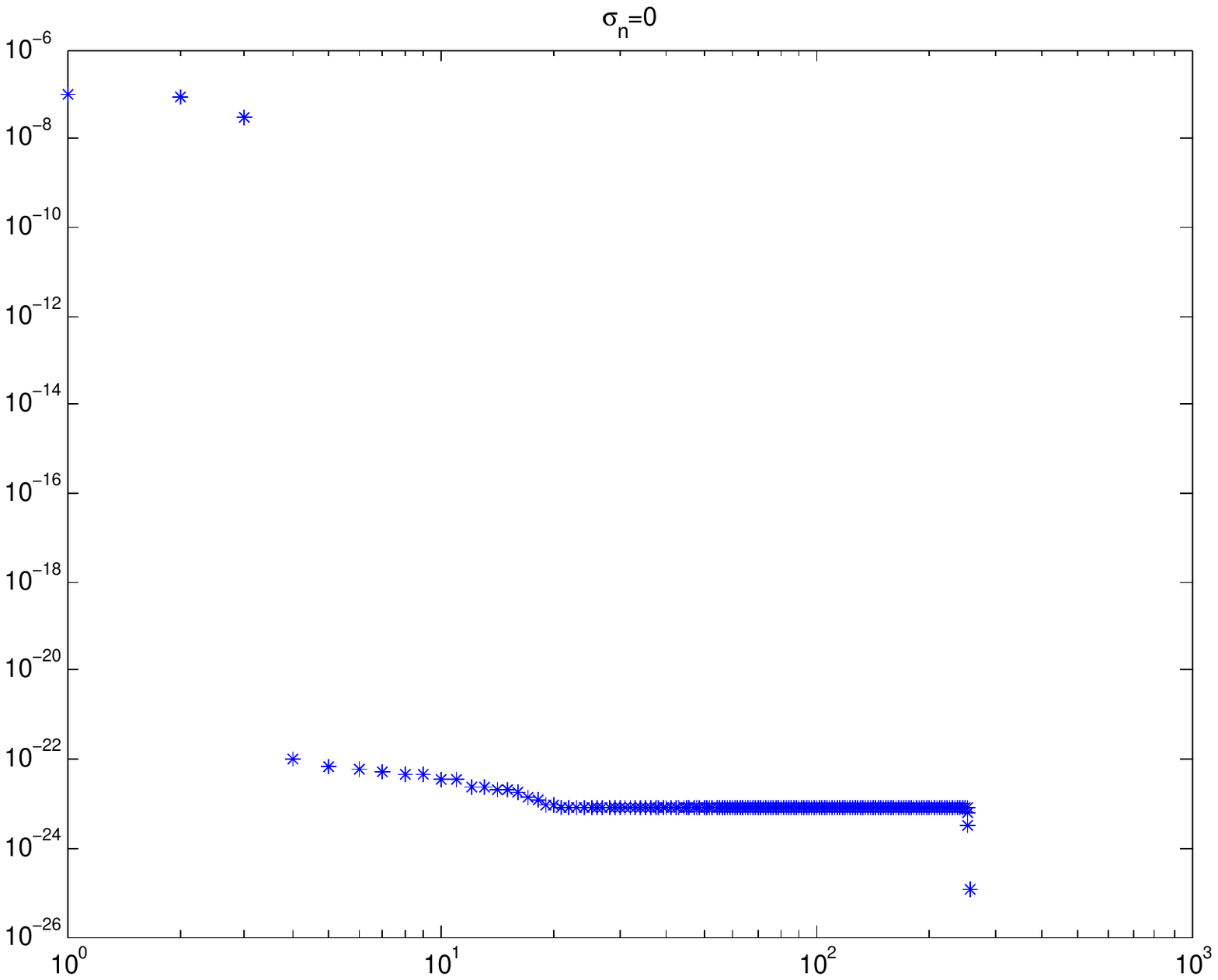}
\includegraphics[width=0.6\textwidth]{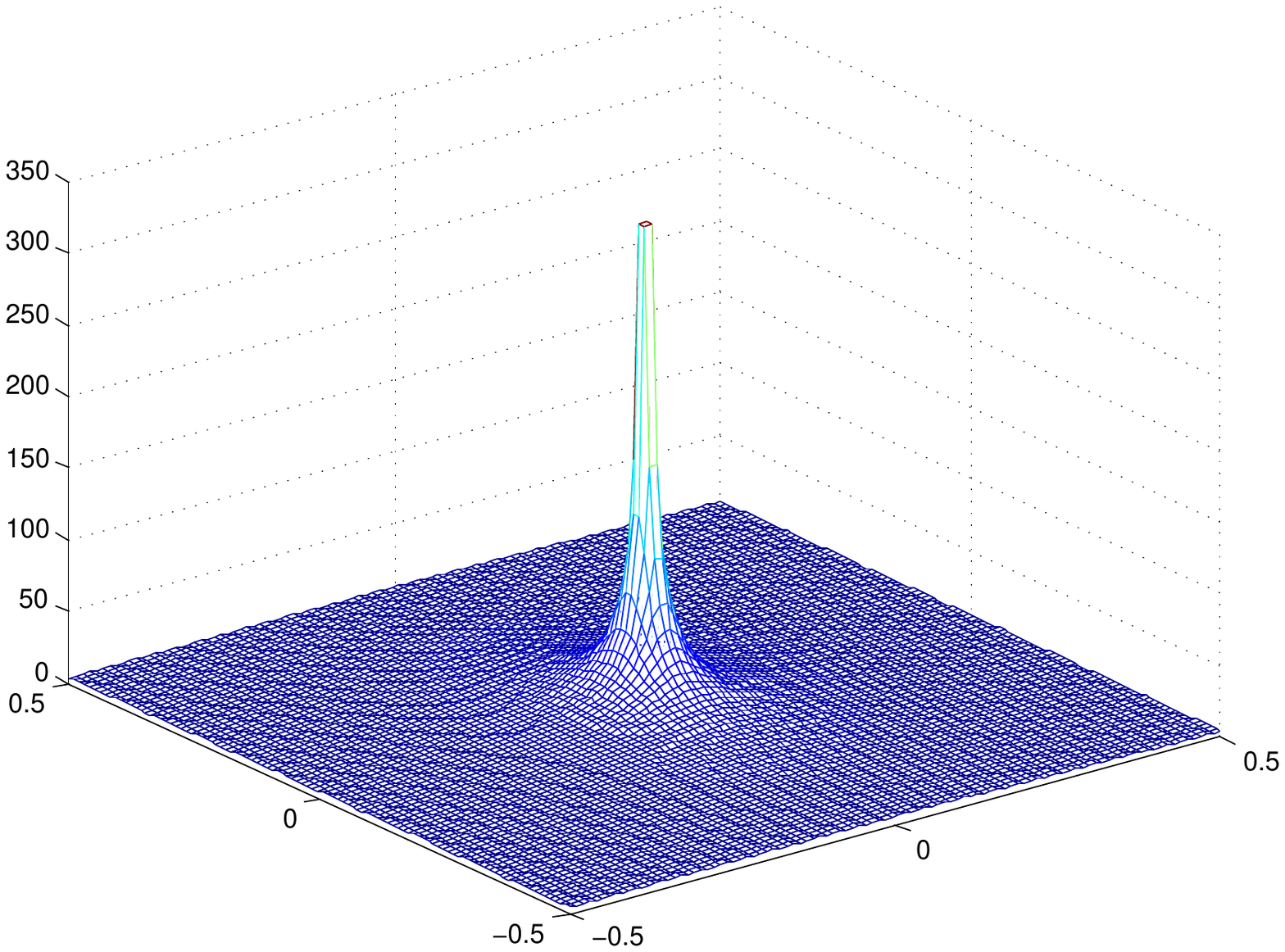}
\caption{\color{black}
Finding the location of an an ellipsoidal target.
The ellipsoid is defined by the equation $x^2+y^2+\frac{z^2}{4}\leq \alpha^2$, $\alpha=0.01$ m.
The electromagnetic parameters are
$\sigma=5.97e7$ S/m, $\mu_*=\mu_0= 1.2566e-06$ H/m, $\omega=133.5$.
The  sources and the receivers are laid out
as indicated by Figure \ref{geometry}.
 Left: log-log plot of the singular values of MSR matrix defined in (\ref{MSR_0}) computed
using approximation formula (\ref{MSR_def}).
Right: magnitude of imaging functional $\mathcal{I}_{MU}$ defined in (\ref{music})  plotted on the $xy$ plane.  
As expected, this functional peaks at  the center of the target.
}\label{localization}
\end{figure}

Figure \ref{localization} shows the localization results. 
The MSR has three dominant singular values indicating that there is only  one target.
The functional $\mathcal{I}_{MU}$ defined in (\ref{music})   peaks at  the center of the target, as anticipated.\\
\color{black}
Next, we  assess  
how this MUSIC location algorithm is capable of differentiating two distinct targets.
%
To do that, pick two small targets 
shaped as previously and centered in the $xy$ plane.
Denote by $\bz_1$ and $\bz_2$ their centers, assume that
$\bz_1$ is at the origin (accordingly the second  ellipsoid
is given by the equation $
(x-   \bz_2\cdot \bm{e}_1 )^2+(y-   \bz_2\cdot \bm{e}_2 )^2+\frac{z^2}{4}\leq \alpha^2$, $\alpha=0.01$ m).  
These two targets have same conductivity $\sigma=5.97e7$ S/m,
permeability is set to be constant everywhere $\mu_0= 1.2566e-06$ H/m,
and as previously $\omega=133.5$.
In this simulation set $L$ to be a positive distance and place 
(on a uniform grid) 256 sources on the square  $[-2,2]^2\times\{L\}$
256 receivers on the square  $[-2,2]^2\times\{-L\}$.
Denote by $\sigma_1$ the maximum  singular value of the MSR matrix without noise, that is, $U\mathcal{M}_q V_p$.
We define the signal-to-noise ratio by 
$$\mathrm{SNR} =\sigma_1/\sigma_{\mathrm{noise}},$$
and the noise level as $\mathrm{SNR}^{-1}$. 

\begin{figure}
\includegraphics[width=\textwidth]{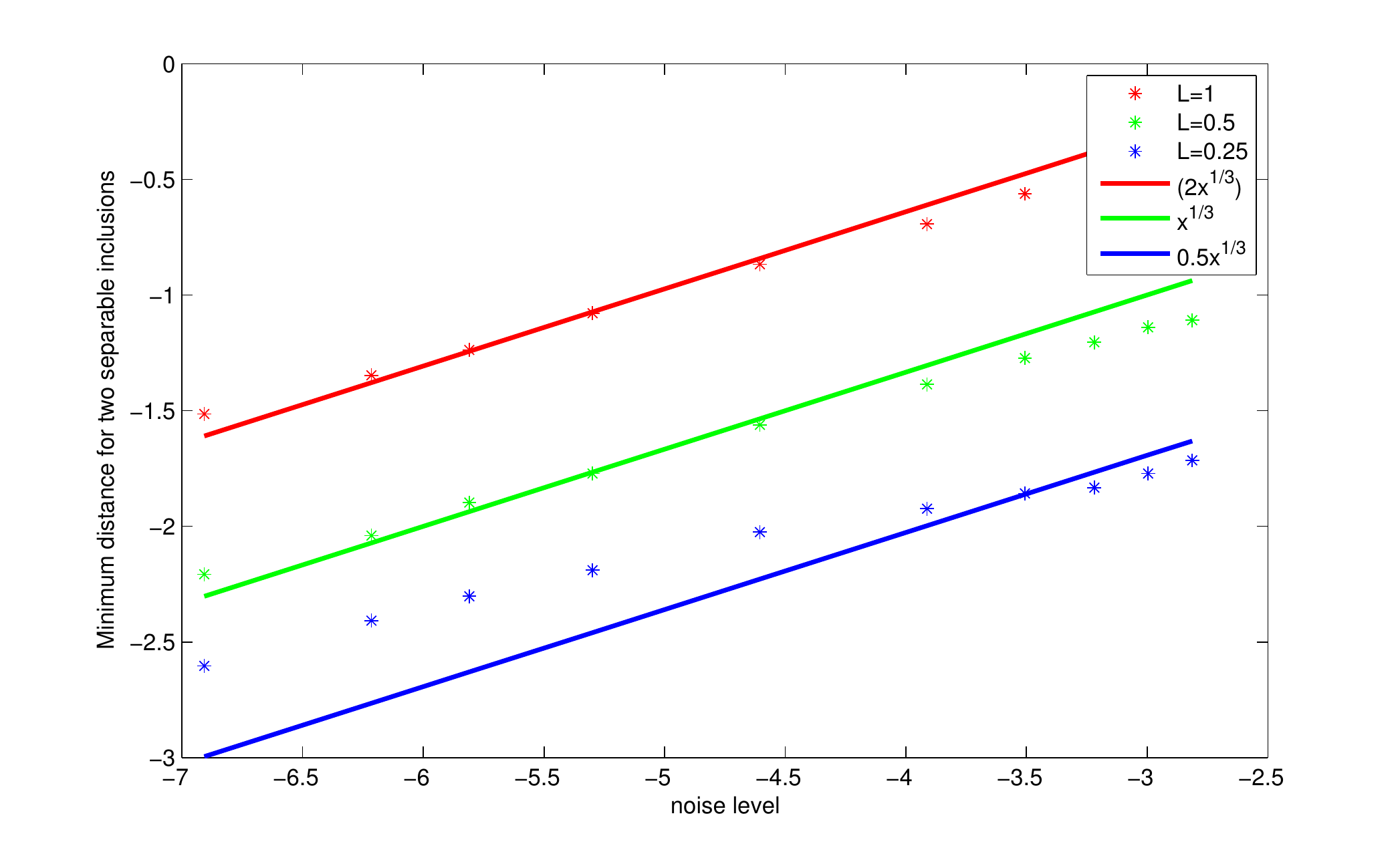}
\caption{
${d_{\mathrm{min}}}$, the minimum distance  needed to 
 differentiate the two ellipsoidal targets with diameter 0.02, plotted as stars against the $\mathrm{SNR}$. 
Units are logarithmic. 
The solid lines illustrate the estimate 
${d_{\mathrm{min}}} \sim 2L \; \mathrm{SNR}^{-1/3}$. 
}\label{fig-res}
\end{figure}

In Tables \ref{tab-res1}-\ref{tab-res5}, we give  for different values of $L$ the minimum distance $\textcolor{black}{d_{\mathrm{min}}}$ \textcolor{black}{between $\bz_1$ and $\bz_2$} needed to 
clearly differentiate the two targets. In Figure \ref{fig-res}, we plot the minimum  distance $\textcolor{black}{d_{\mathrm{min}}}$ against
 $\mathrm{SNR}$ for $L=1.0,0.5, 0.25$ in logarithmic scale. 
We observe that the minimum distance $\textcolor{black}{d_{\mathrm{min}}}$ is approximately  equal to $2L \; \mathrm{SNR}^{-1/3}$.

\begin{table}
\begin{tabular}{|c|c|c|c|c|c|c|c|c|c|c|}
\hline
noise level&0.1\%&0.2\%&0.3\%&0.5\%& 1\% & 2\% & 3\% & 4\% &5\% & 6\%\\
\hline
$\textcolor{black}{d_{\mathrm{min}}}$&0.27&0.33&0.37&0.43 & 0.53 & 0.63 & 0.68 &0.74 & 0.78& 0.84\\
\hline
\end{tabular}
\caption{Computed values of ${d_{\mathrm{min}}}$, the minimum distance  needed to 
 differentiate the two ellipsoidal targets with diameter $0.02$.
Here $L=1.25$, where $L$ is the distance from the targets to the plane containing the sources,
which we chose to be equal to the distance from the targets to the plane containing the receivers.
}\label{tab-res1}
\end{table}

\begin{table}
\begin{tabular}{|c|c|c|c|c|c|c|c|c|c|c|}
\hline
noise level &0.1\%&0.2\%&0.3\%&0.5\%& 1\% & 2\% & 3\% & 4\% &5\% & 6\%\\
\hline
$\textcolor{black}{d_{\mathrm{min}}}$&0.22&0.26&0.29&0.34 & 0.42 & 0.50 & 0.57 &0.62 & 0.66& 0.69\\
\hline
\end{tabular}
\caption{Same as previous table  for $L=1$. }\label{tab-res2}
\end{table}

\begin{table}
\begin{tabular}{|c|c|c|c|c|c|c|c|c|c|c|}
\hline
noise level&0.1\%&0.2\%&0.3\%&0.5\%& 1\% & 2\% & 3\% & 4\% &5\% & 6\%\\
\hline
$\textcolor{black}{d_{\mathrm{min}}}$&0.16&0.19&0.22&0.26 & 0.31 & 0.38 & 0.42 &0.44 & 0.48& 0.50\\
\hline
\end{tabular}
\caption{Same as previous table for $L=0.75$. }\label{tab-res3}
\end{table}

\begin{table}
\begin{tabular}{|c|c|c|c|c|c|c|c|c|c|c|}
\hline
noise level&0.1\%&0.2\%&0.3\%&0.5\%& 1\% & 2\% & 3\% & 4\% &5\% & 6\%\\
\hline
$\textcolor{black}{d_{\mathrm{min}}}$&0.11&0.13&0.15&0.17 & 0.21 & 0.25 & 0.28 &0.30 & 0.32& 0.33\\
\hline
\end{tabular}
\caption{Same as previous table for $L=0.5$. }\label{tab-res4}
\end{table}

\begin{table}
\begin{tabular}{|c|c|c|c|c|c|c|c|c|c|c|}
\hline
 noise level&0.1\%&0.2\%&0.3\%&0.5\%& 1\% & 2\% & 3\% & 4\% &5\% & 6\%\\
\hline
$\textcolor{black}{d_{\mathrm{min}}}$&0.074&0.09&0.1&0.112 & 0.132 & 0.146 & 0.156 &0.16 & 0.17& 0.18\\
\hline
\end{tabular}
\caption{Same as previous table for $L=0.25$. }\label{tab-res5}
\end{table}

\subsection{Performance of the classification algorithm}

\begin{figure}
\includegraphics[width=0.16\textwidth]{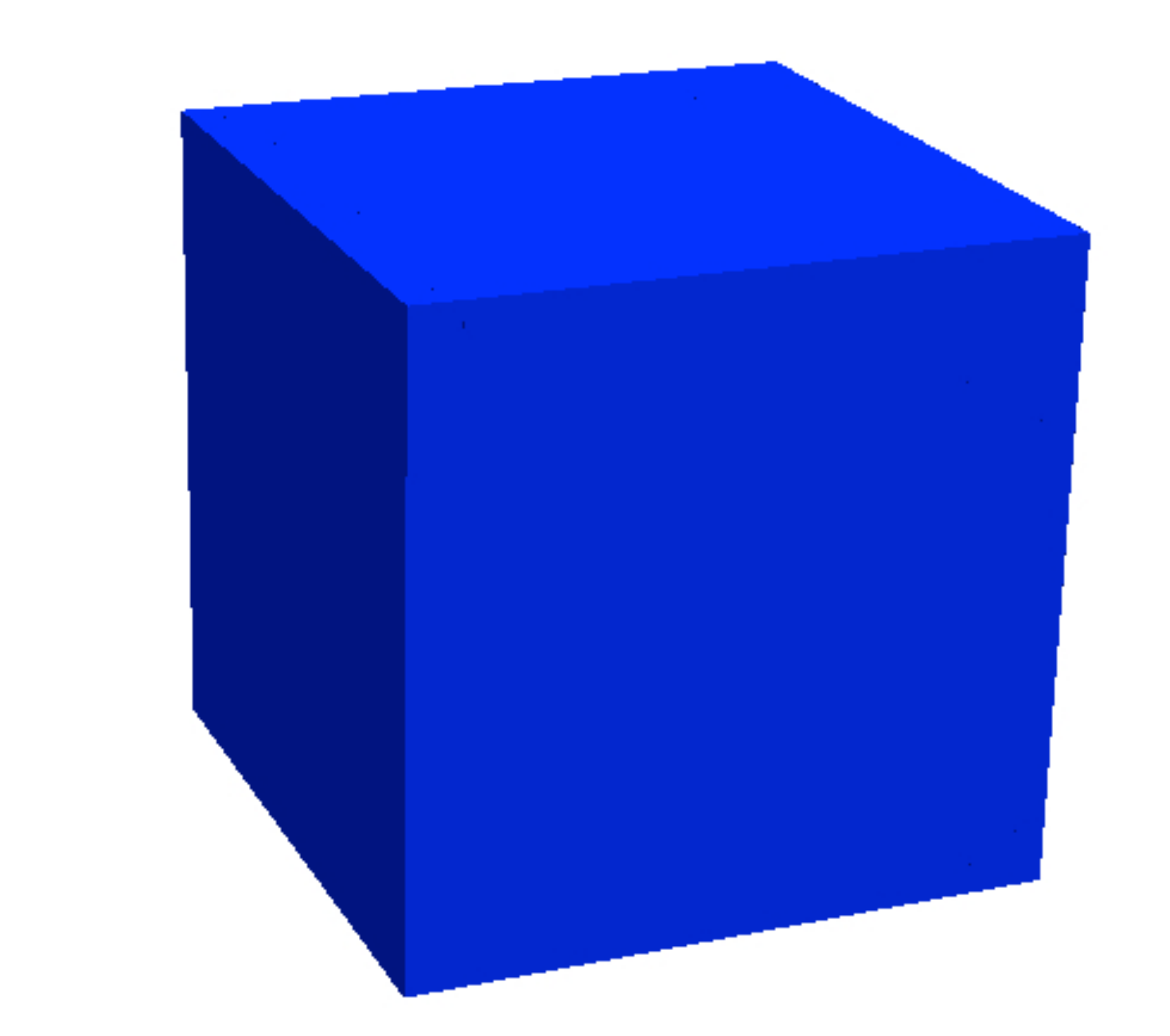}
\includegraphics[width=0.16\textwidth]{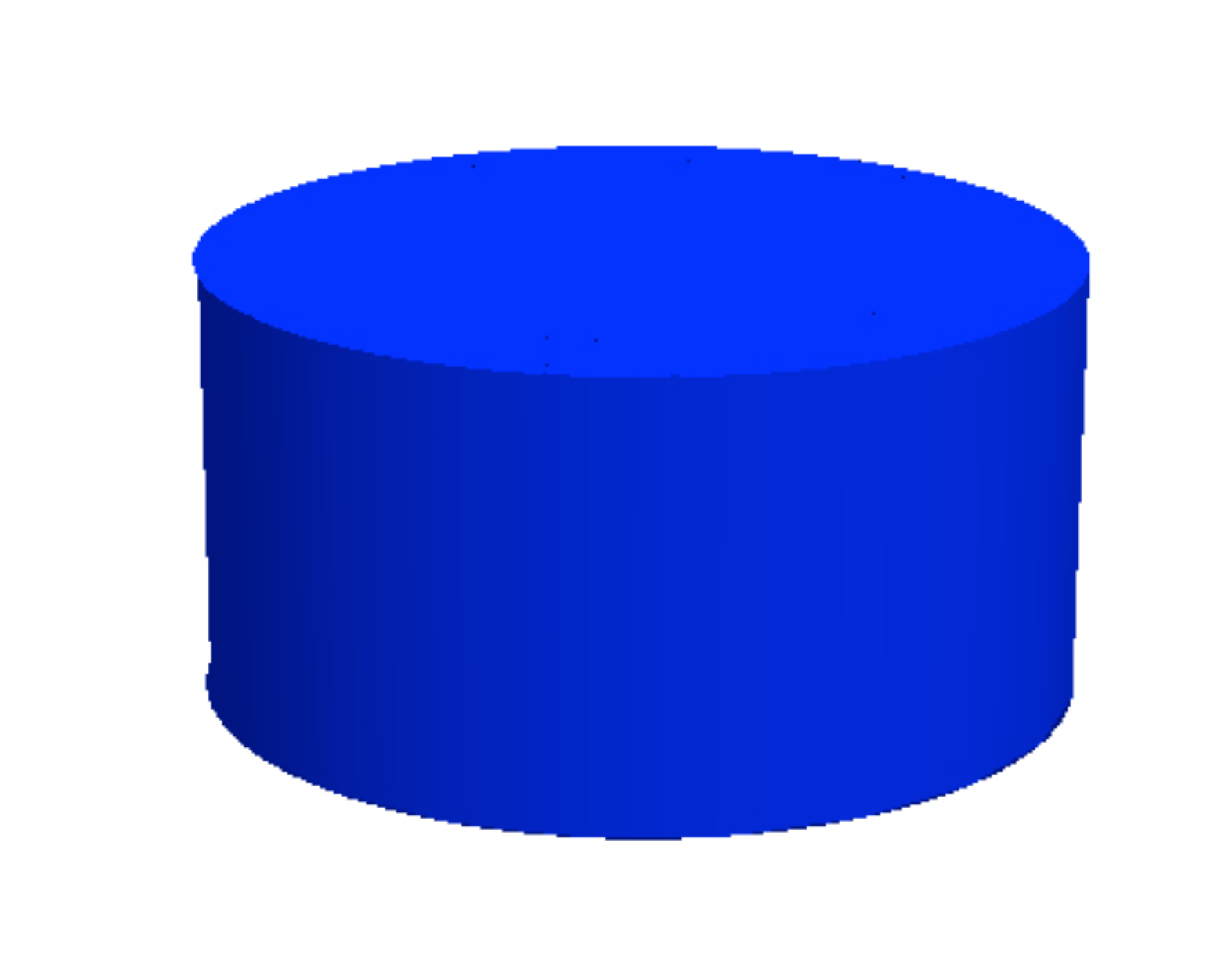}
\includegraphics[width=0.16\textwidth]{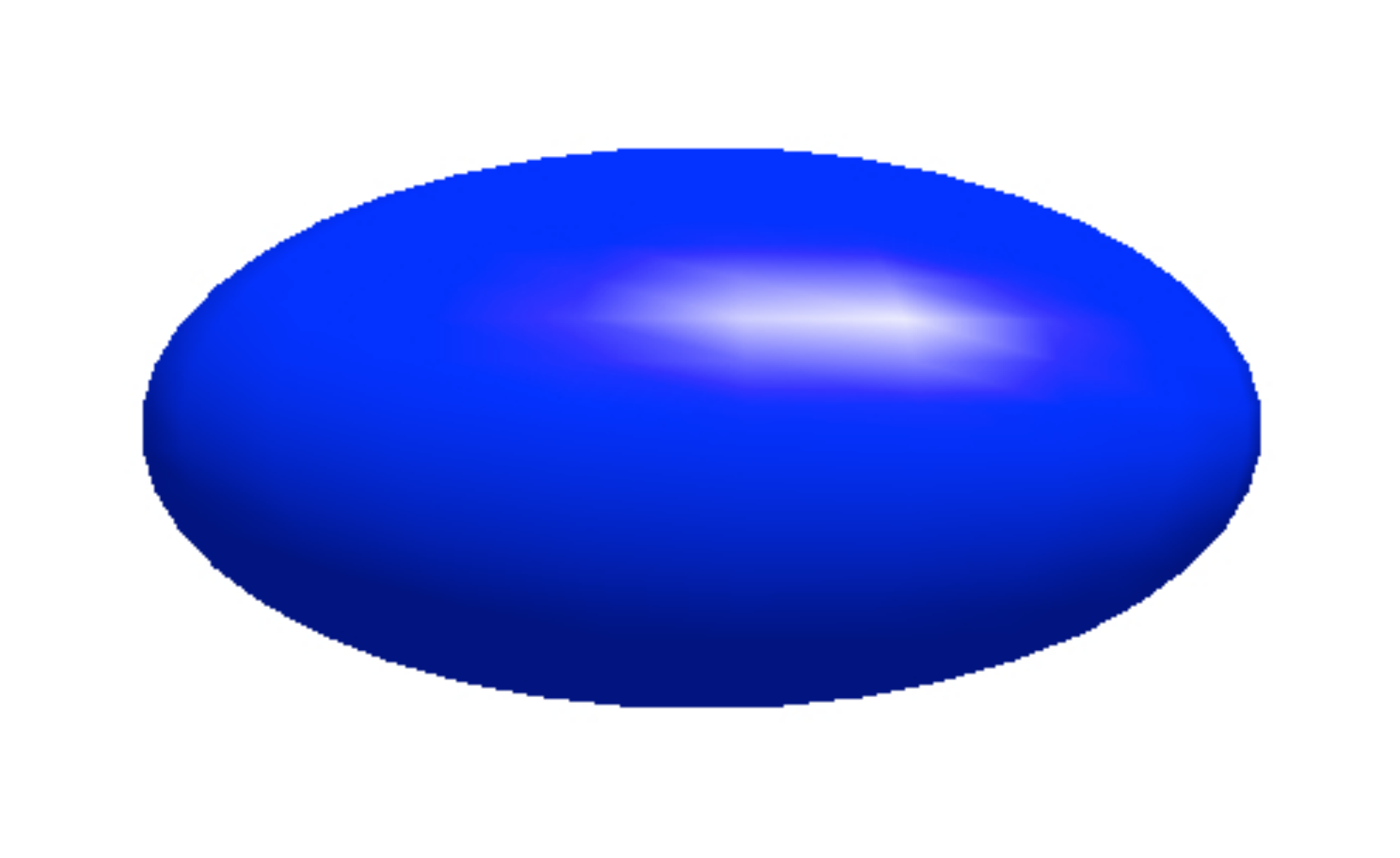}
\includegraphics[width=0.16\textwidth]{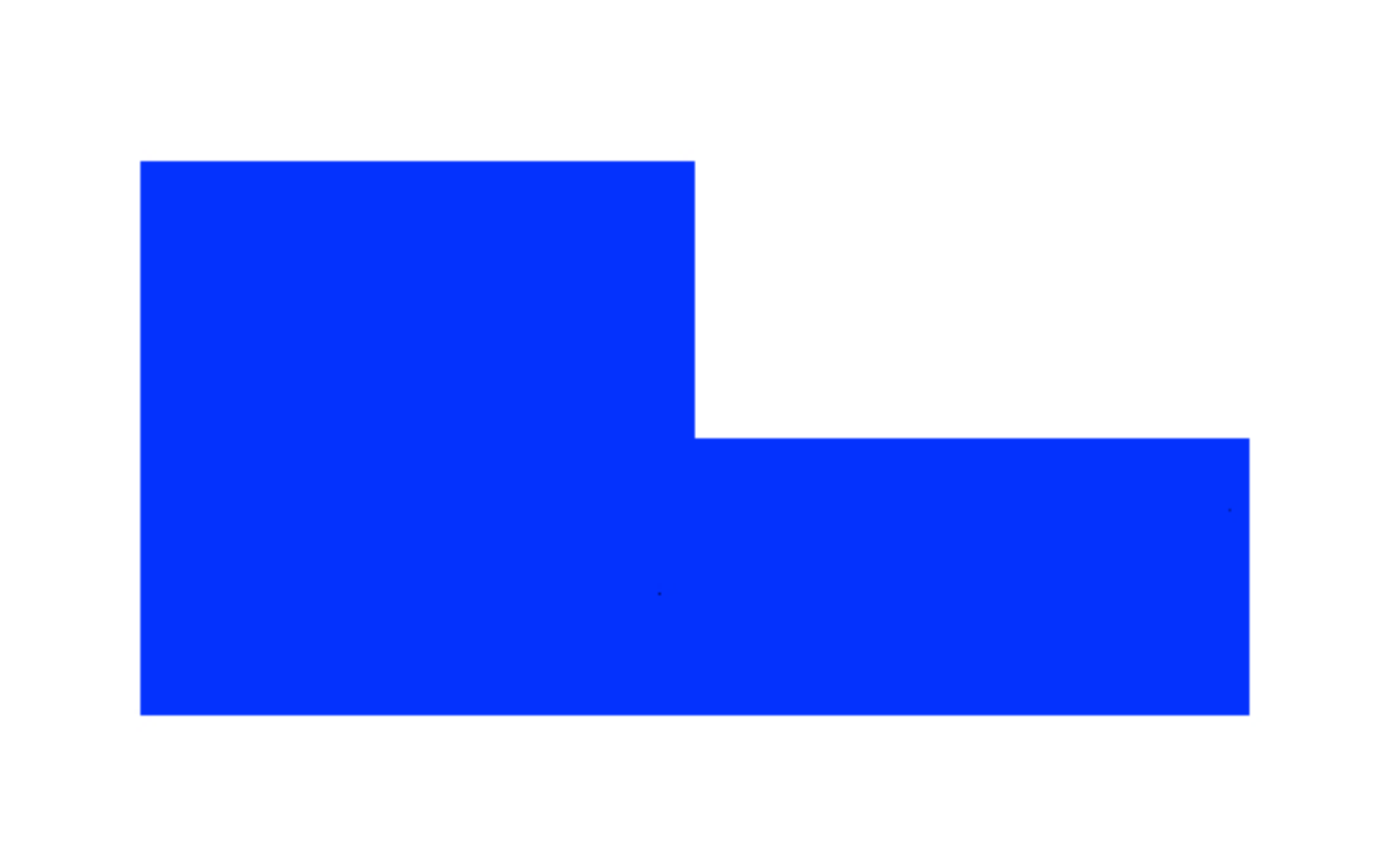}
\includegraphics[width=0.16\textwidth]{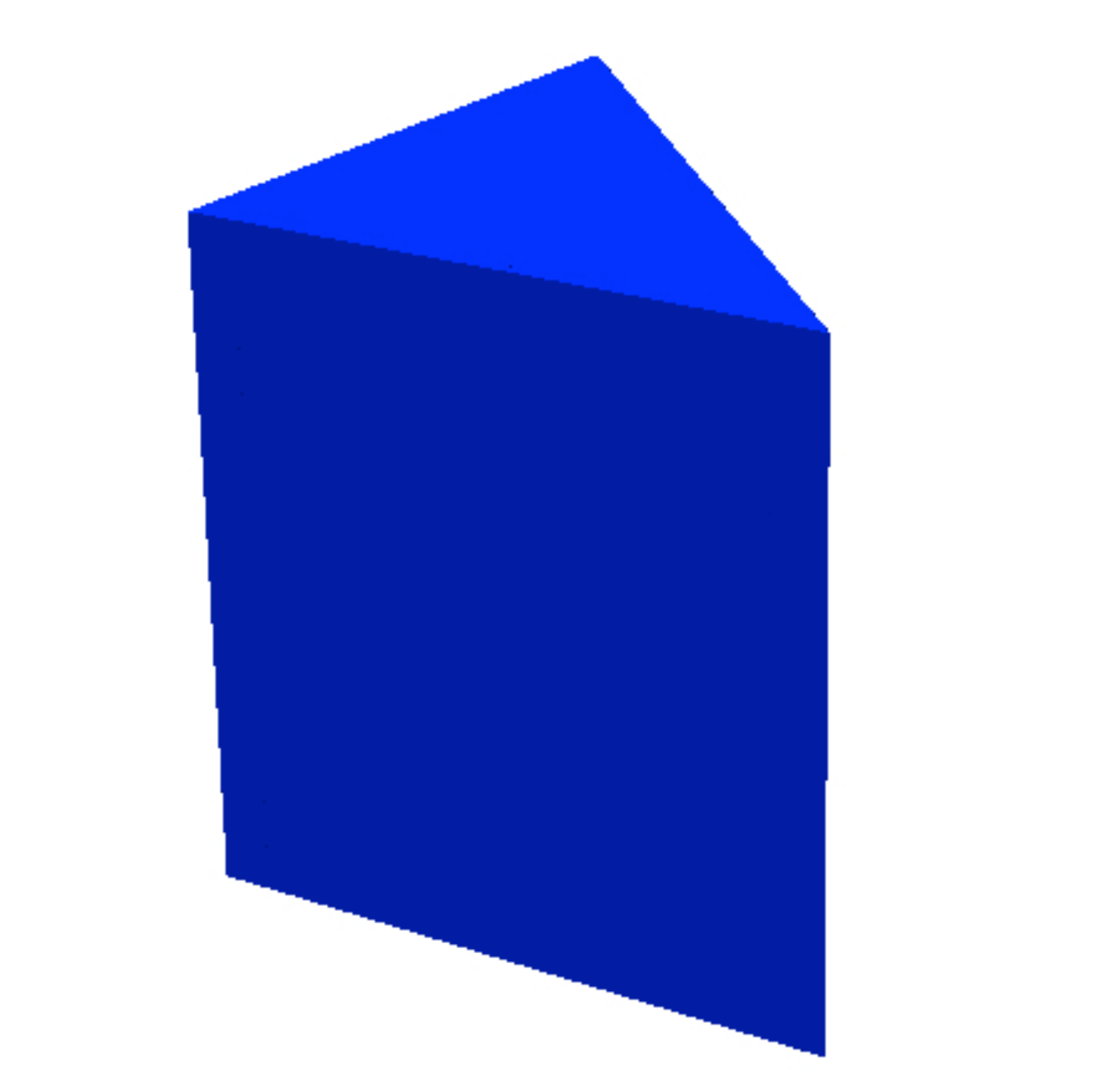}
\includegraphics[width=0.16\textwidth]{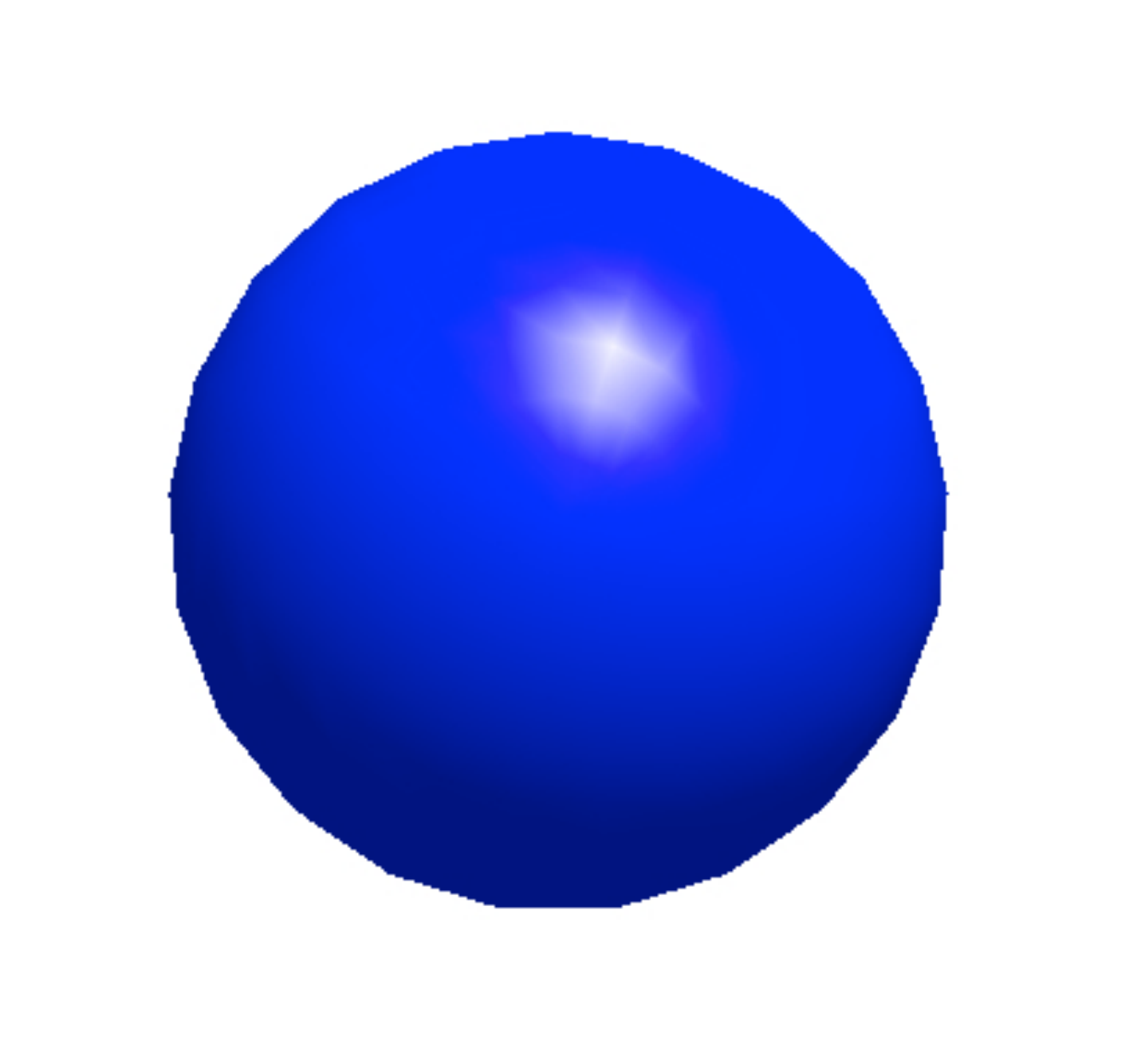}
\caption{The shapes in the dictionary.}\label{fig-dict}
\end{figure}

Next, we  report some numerical results to demonstrate the efficiency of Algorithm \ref{alg-1} at a single
frequency and at multiple frequencies.  The dictionary includes the following domains: (1) cube
$[-1,1]^3$, (2) cylinder $\{x^2+y^2\leq 1, -0.5\leq z\leq 0.5\}$, (3) ellipsoid $x^2+y^2+\frac{z^2}{4}\leq1$, (4) L-shaped domain $[-1,1]\times[-0.5,0.5]\times[-0.5,0.5]$, (5) prism $\{-1\leq x, -1\leq y, \mbox{ and } x+y\leq 1\}\times [-1,1]$, and (6) sphere $x^2+y^2+z^2\leq1$. These shapes are sketched in Figure \ref{fig-dict}.

\subsection{Classification from measurements at a single frequency}
Table \ref{tab-svd} indicates
for each of the five domains listed above
 the three significant  singular values of $\mathbb{M}$  at the operating frequency $\omega=133.5$.
\begin{table}
\centering{
\begin{tabular}{|c|c|}
\hline
shape & singular values\\
\hline
cube & 2.2485, 2.2485, 2.2484\\
\hline
cylinder & 0.5997, 0.5997,0.3429\\
\hline
ellipsoid &2.6159,2.1916,2.1916\\
\hline
L-shape &0.1316,0.1278,0.0941\\
\hline
prism    & 3.0423, 2.8299,2.3296\\
\hline
sphere  & 0.8282, 0.8277,0.8277\\
\hline
\end{tabular}
\begin{tabular}{|c|c|}
\hline
shape & singular values\\
\hline
cube & 1.0, 1.0, 1.0\\
\hline
cylinder & 1.0, 1.0,0.5717\\
\hline
ellipsoid &1.0,0.8378,0.8377\\
\hline
L-shape &1.0,0.9715,0.7151\\
\hline
prism    & 1.0, 0.9302,0.7657\\
\hline
sphere  & 1.0, 0.9993,0.9993\\
\hline
\end{tabular}
}
\caption{Three significant singular values of $\mathbb{M}$ for targets with different shapes. The  table on the left is for the original  singular values while the table on the right is for the normalized singular values.}\label{tab-svd}
\end{table}

\begin{figure}
\includegraphics[scale=.3]{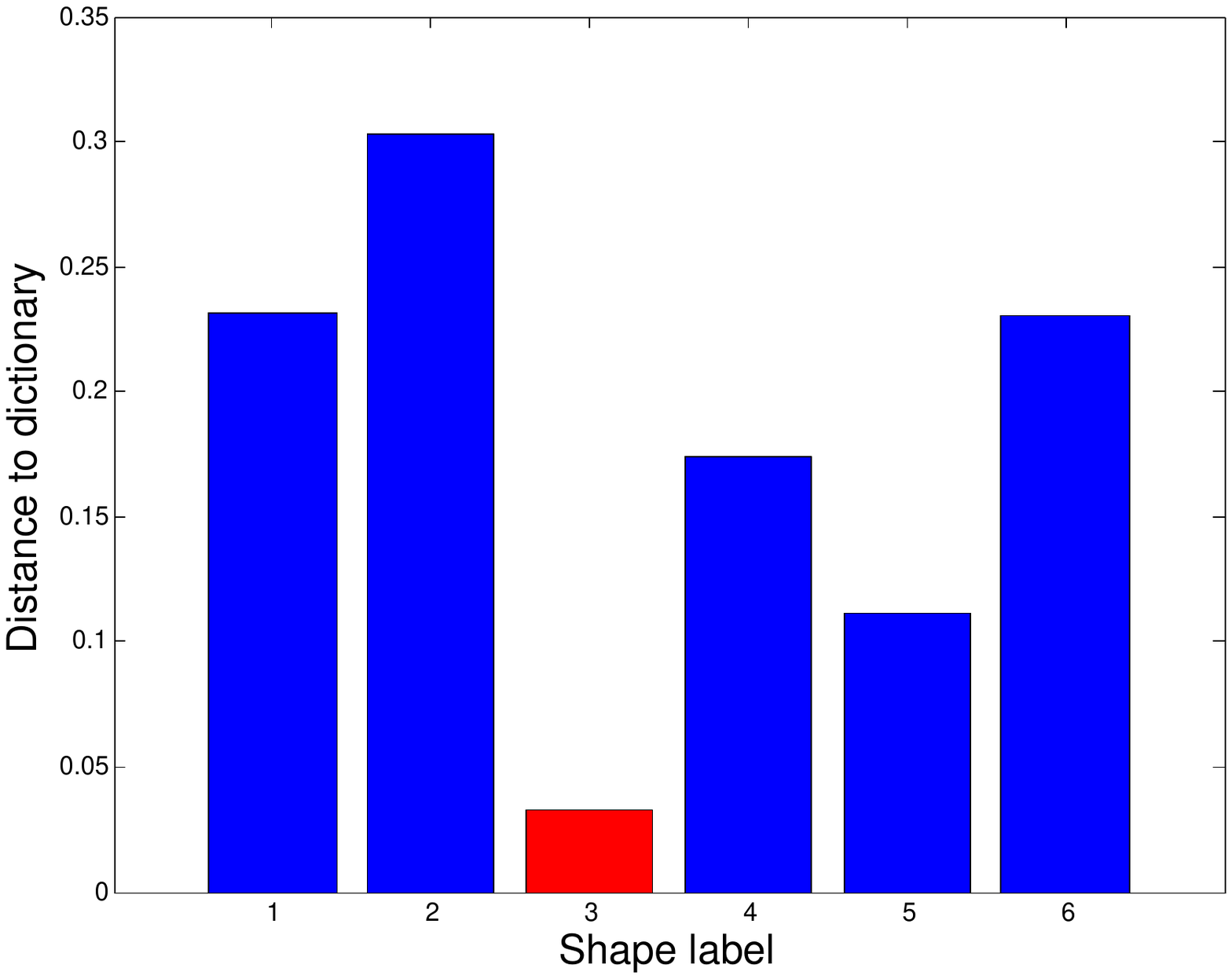}
\includegraphics[scale=.3]{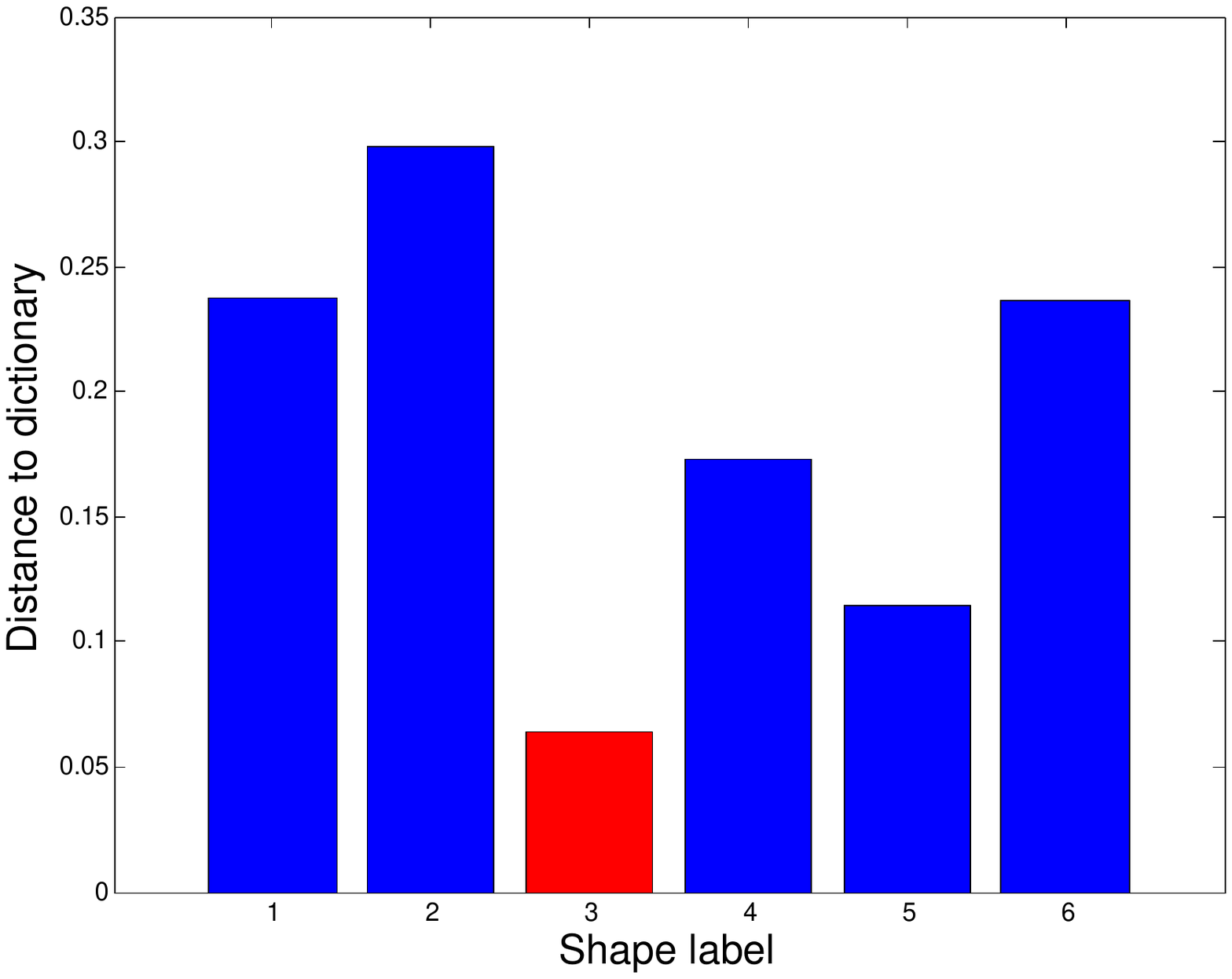}
\includegraphics[scale=.3]{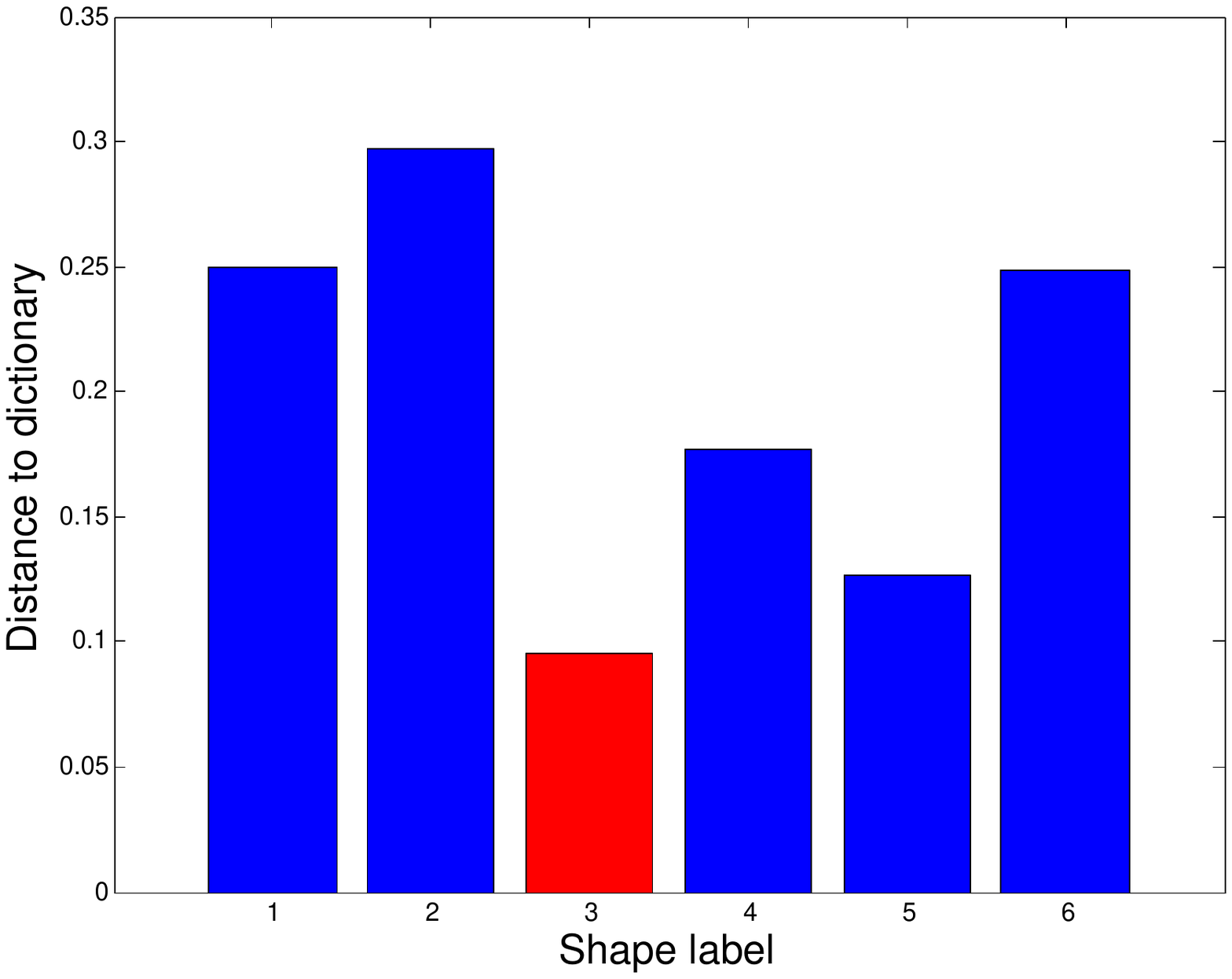}
\includegraphics[scale=.3]{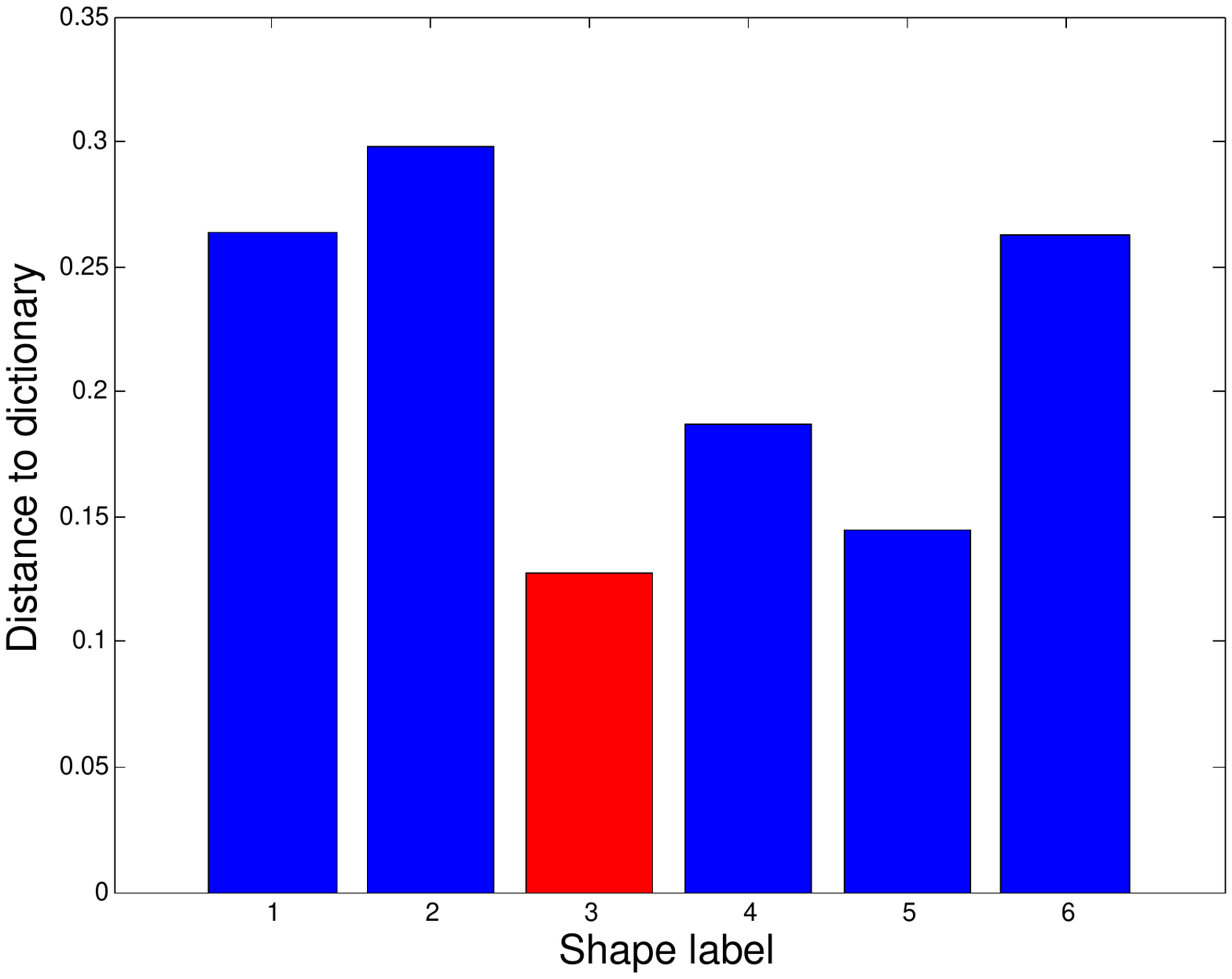}
\caption{Matching results for a rotated ellipsoid target, from left to right, top to bottom, $nl=10\%,20\%,30\%,40\%$. 
\color{black} Labels on horizontal axis: 
1 -cube, 2 -cylinder, 3 -ellipsoid, 4 -L-shaped domain, 5 -prism, 6 -sphere.
Vertical axis: \color{black}distance between the shape descriptors of the target computed from the measurements and the shape descriptors of dictionary.  
The distance is computed by averaging 1000 realizations  for each noise level.}\label{fig-dict-match-rot}
\end{figure}

We first show  results on classification using a single frequency  ($\omega= 133.5$).  This frequency satisfies  $k\alpha^2=O(1)$
so the asymptotic formula in Section \ref{sec2} can be safely used. 
We first  locate the target  by applying the MUSIC algorithm. We then  recover the CPT matrices by 
solving the minimization problem \eqref{lsq-msr}.  
Set $\mathcal{D}=\{D_1, D_2, D_3, D_4, D_5,D_6\}$, 
where these shapes correspond  to the aforementioned five domains labeled in the same order. 
Assume that the number of sources and receivers are both 256. Accordingly, the MSR matrix  is $256$ by $256$. 
In this simulation \color{black} we place these sources on a uniform grid on the square 
$[-2,2]^2\times\{1\}$ and these receivers on a uniform grid on the square 
$[-2,2]^2\times\{-1\}$:
 \color{black} see in Figure \ref{geometry} a sketch  of the geometry for numerical simulations of 
shape detection and classification in this paper.
\color{black}
\color{black} Figure \ref{fig-dict-match-rot} shows the matching results for a target 
whose shape is defined by the equation
 $x^2/4+y^2+z^2\leq \alpha^2$. Note that this
  is just a rotation of the ellipsoidal target from the dictionary. 
  
\begin{figure}
\includegraphics[scale=.3]{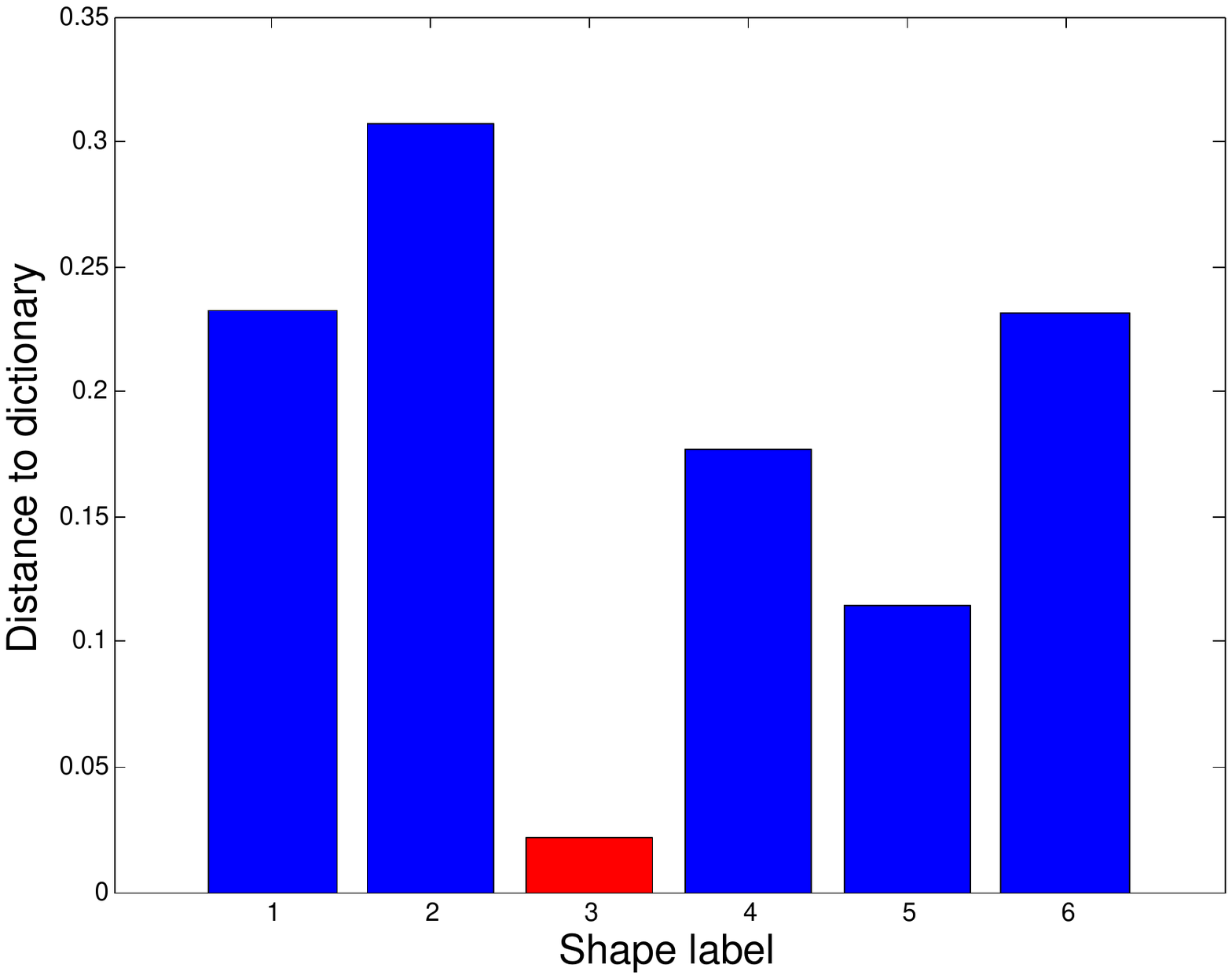}
\includegraphics[scale=.3]{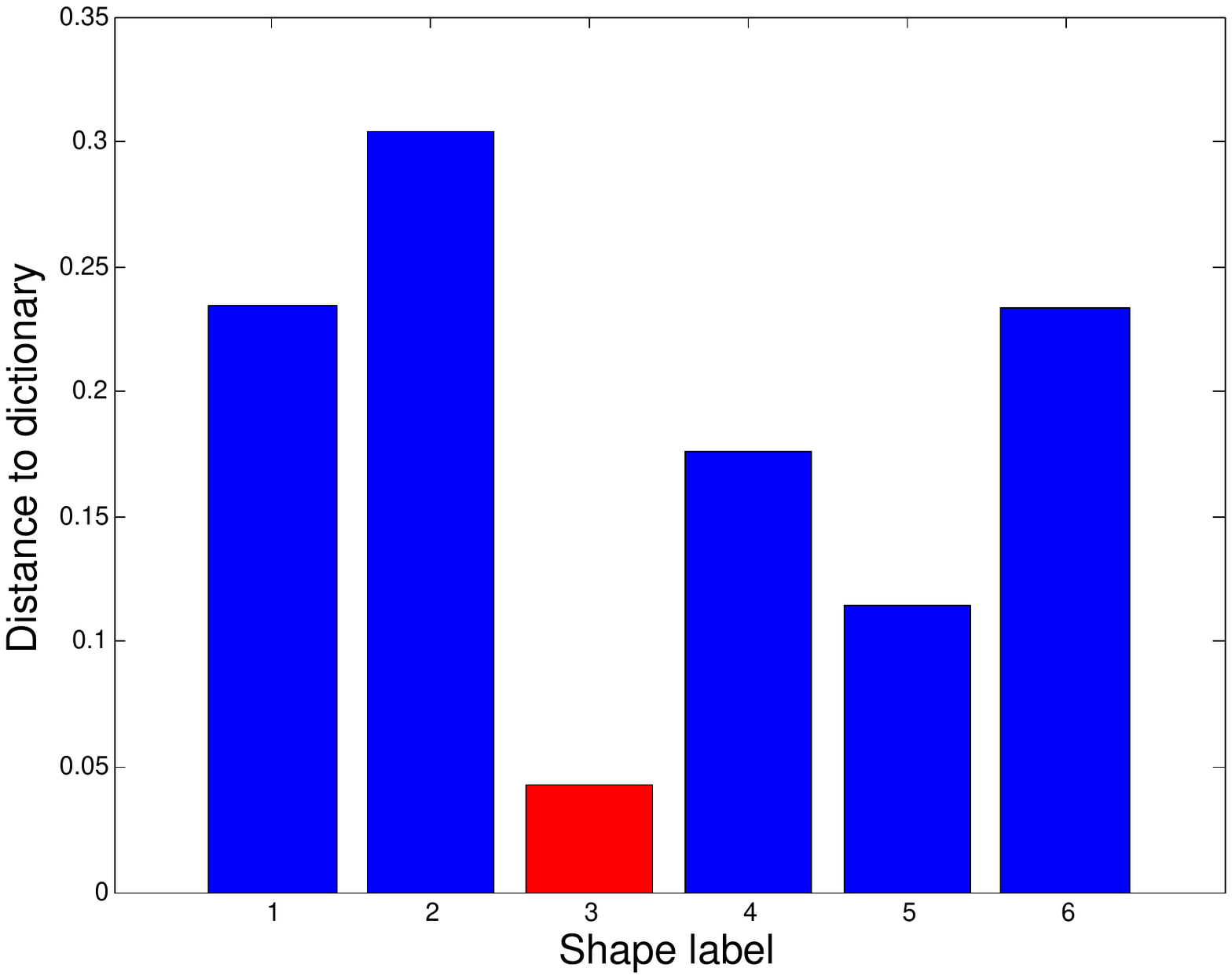}
\includegraphics[scale=.3]{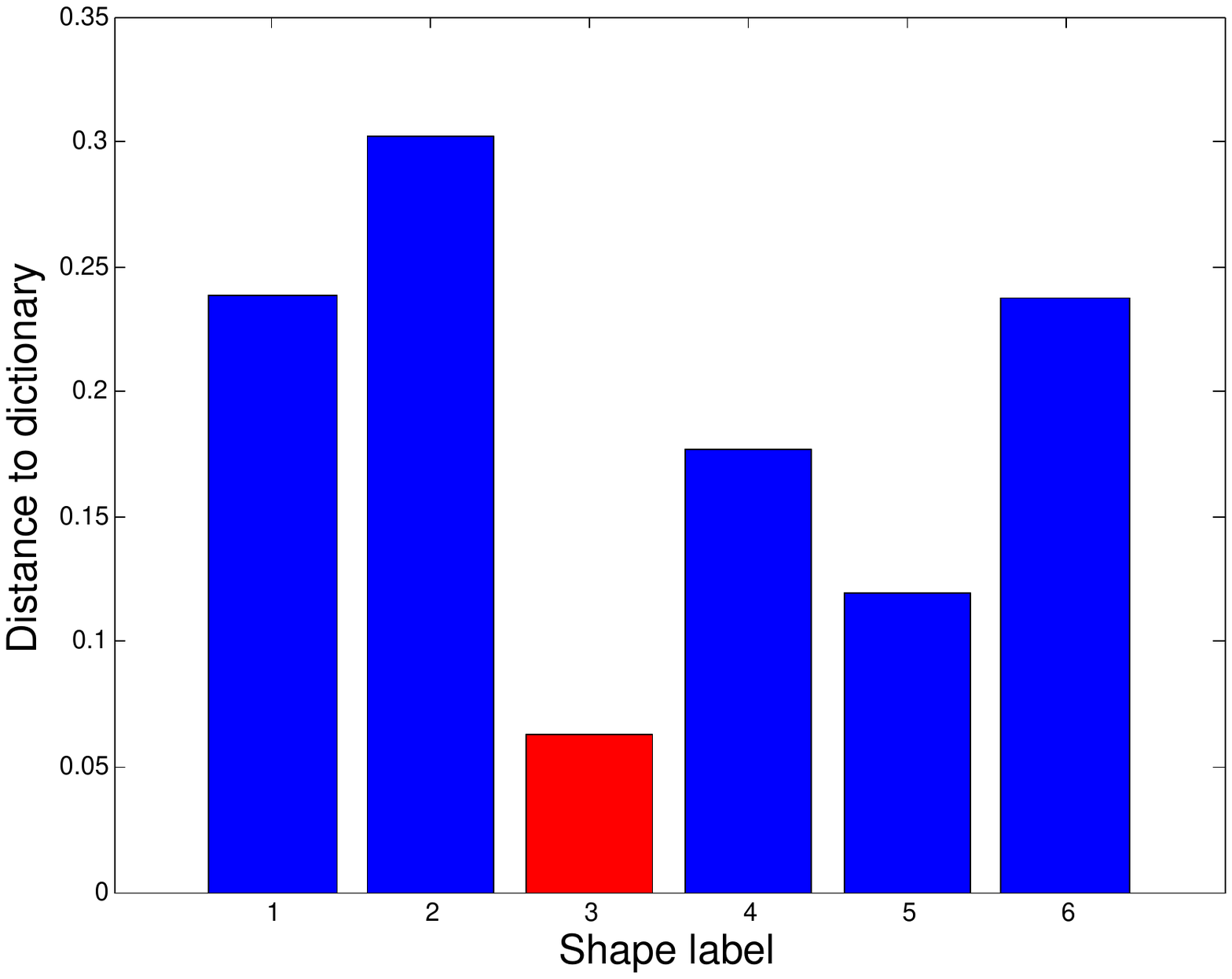}
\includegraphics[scale=.3]{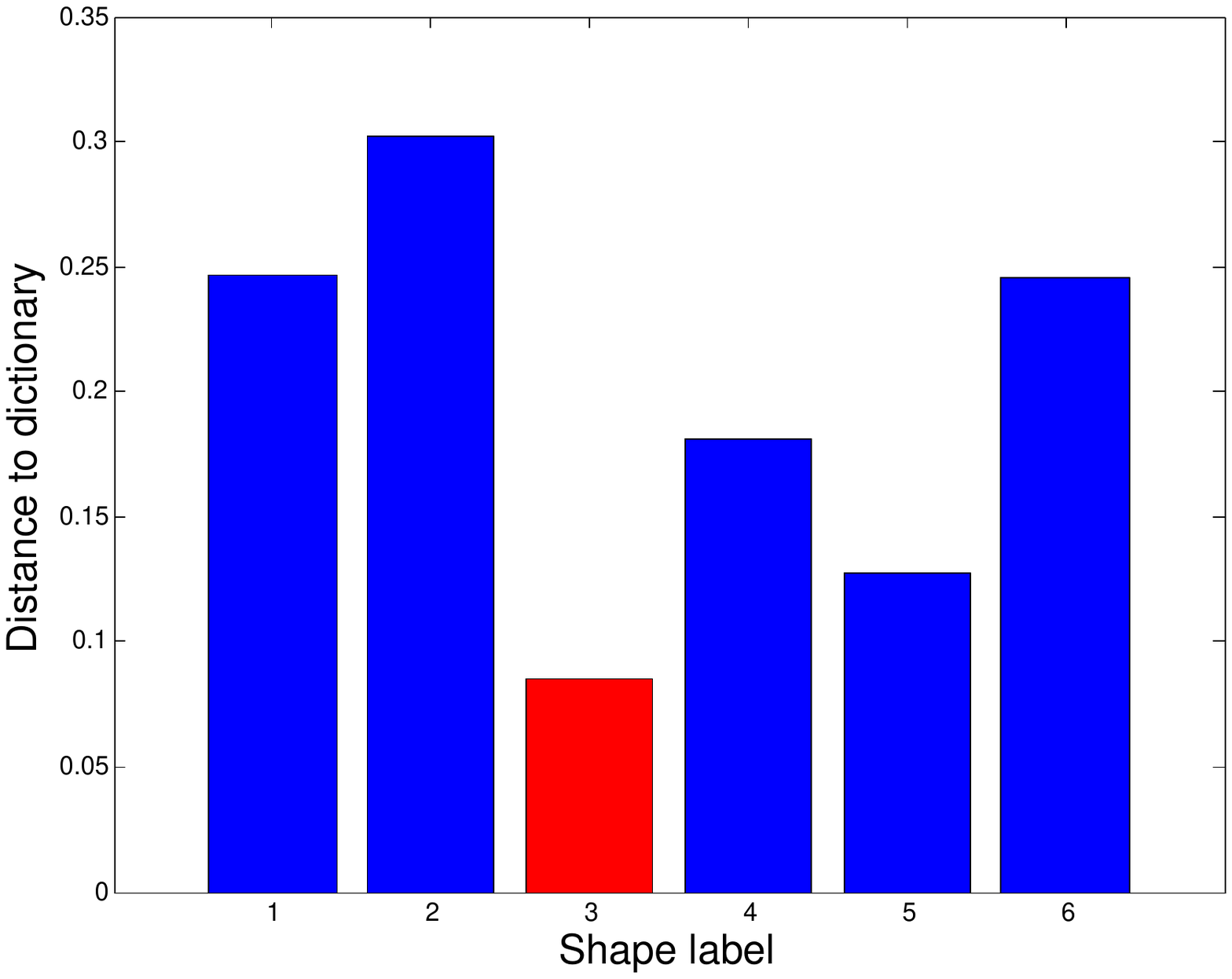}
\caption{Matching results for a small ellipsoid target,  from left to right, top to bottom, $nl=10\%,20\%,30\%,40\%$. 
\color{black} Labels on horizontal axis: 
1 -cube, 2 -cylinder, 3 -ellipsoid, 4 -L-shaped domain, 5 -prism, 6 -sphere.
Vertical axis: \color{black} distance between the shape descriptors of the target computed from the measurements and the shape descriptors of dictionary.  
\color{black}
The distance is computed by averaging 1000 realizations  for each noise level.
}\label{fig-dict-match-small}
\end{figure}

Figure \ref{fig-dict-match-small} shows the results of classification for a small ellipsoid target described by $x^2+y^2+z^2/4\leq 0.25\alpha^2$. \\
 It is visible on this  figure that our algorithm can recognize the correct shape.  
We  find that, at each noise level from our selection,  the minimum 
 $\min_{D\in \mathcal{D}}\{\|D-\hat{D}\|_2\}$
is achieved at $D=D_3$,  that is at  the ellipsoid shaped target.
  It is remarkable that even in the case when the noise level reaches $40\%$, we can still recognize the target.

\color{black}
\subsection{Classification from measurements at multiple frequencies}
\begin{figure}
\includegraphics[scale=.3]{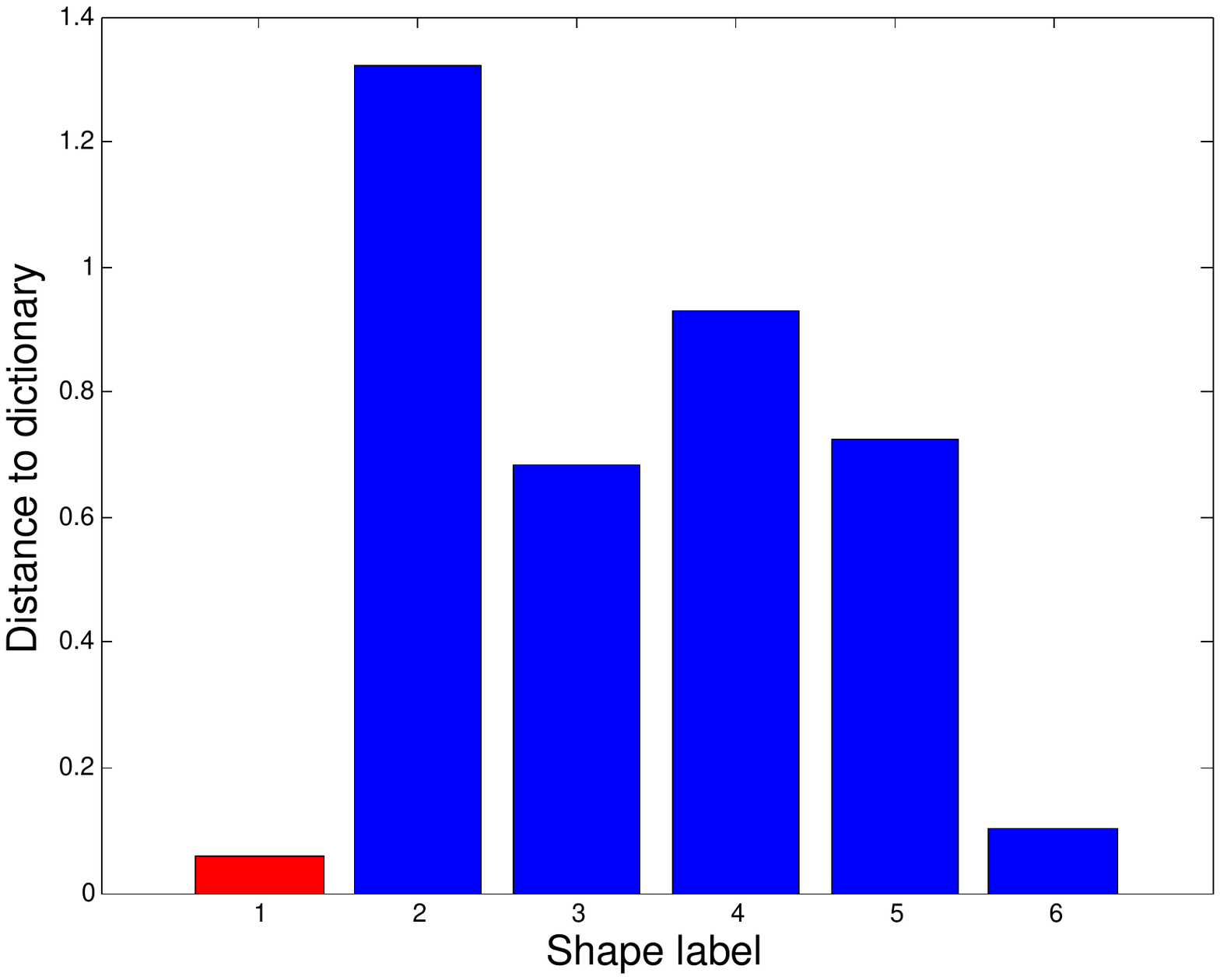}
\includegraphics[scale=.3]{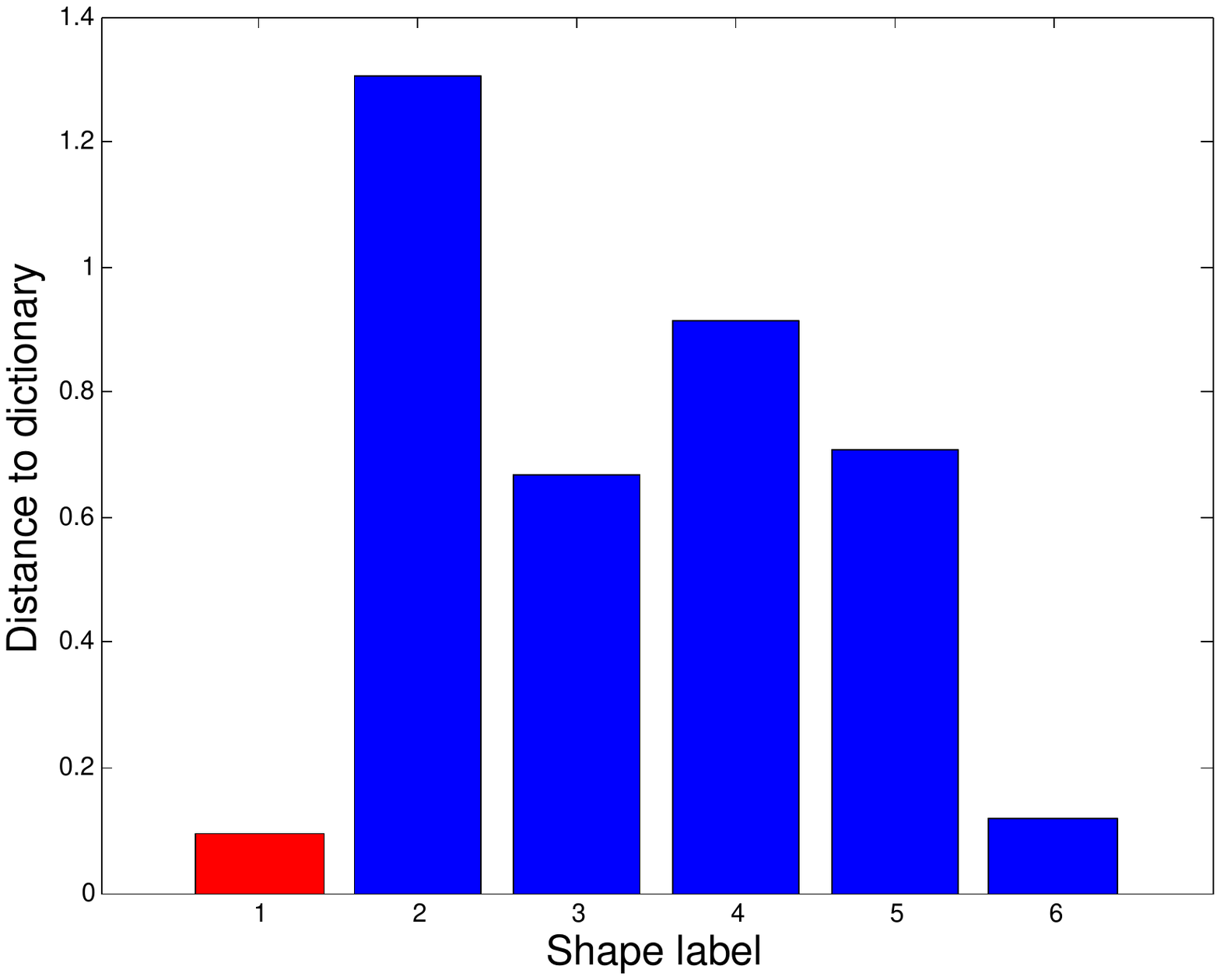}
\includegraphics[scale=.3]{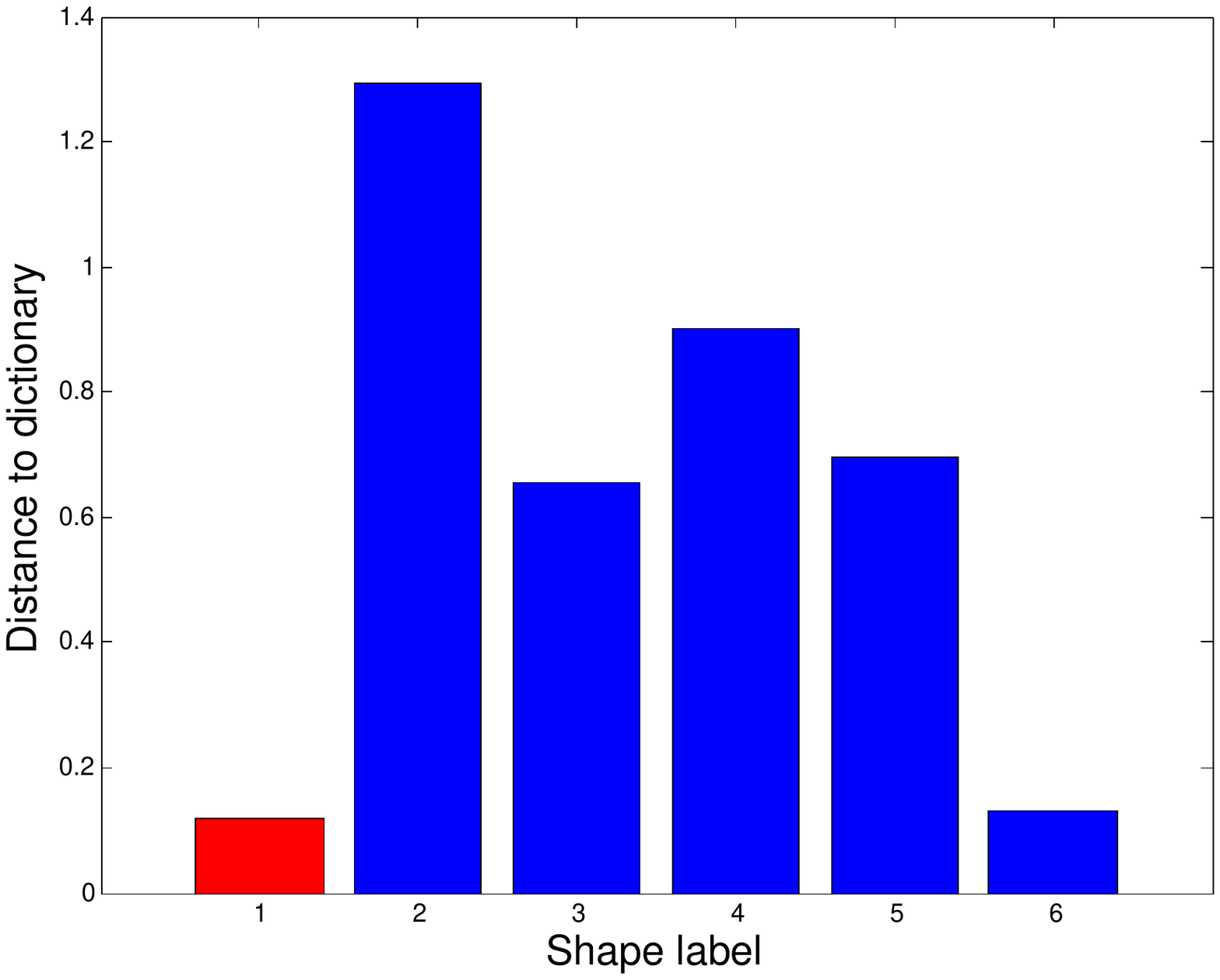}
\includegraphics[scale=.3]{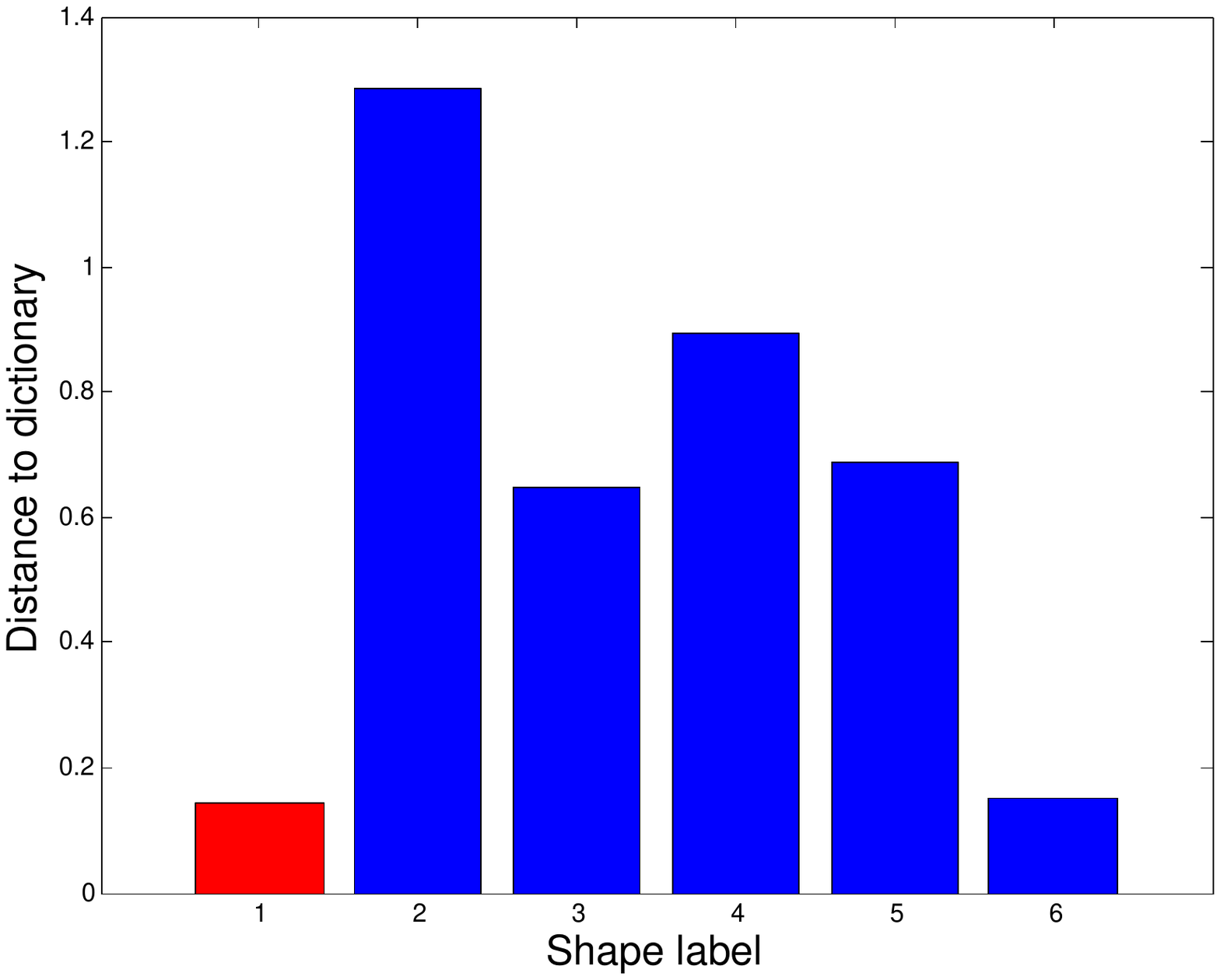}
\caption{\color{black} Use of multiple frequencies for identifying  a cube. 
Frequencies are
$\omega_n=73.5+10n,n=1,2,\ldots,19$, so the requirement $\nu=O(1)$ for 
approximation formula from Theorem \ref{thm-1} to be valid is not violated
(given that we set $\sigma=5.97e7$S/m, $\mu_*=\mu_0= 1.2566e-06$ H/m, and $\alpha=0.01$ m).
Labels on horizontal axis: 
1 -cube, 2 -cylinder, 3 -ellipsoid, 4 -L-shaped domain, 5 -prism, 6 -sphere.
Vertical axis: distance between the shape descriptors of the target computed from the measurements and the shape descriptors of dictionary.  
The distance is computed by averaging 1000 realizations  for each noise level.
From left to right, top to bottom: the noise level is equal to $5\%,8\%,10\%,12\%$.}\label{multi-recog}
\end{figure}

In Table \ref{tab-svd}, we show that the normalized singular values for the cube and the sphere are very similar, so we can not distinguish a cube from a sphere if the data is corrupted by noise. The dependence of CPT matrix $\mathbb{M}$ on $s$ is nonlinear, in other words, $\mathbb{M}$ is nonlinear with respect to $\omega$. This motivates 
 trying to  use multiple frequencies in order to 
be able to 
differentiate them. In our simulation we used the frequencies
$\omega_n=73.5+10n,n=1,2,\ldots,19$. The highest frequency in this range is 
$\omega=263.5$, yielding $k\alpha^2\approx 2$:
the skin depth $\delta$ is close to $ \alpha$ and our basic asymptotic approximation is still valid. 
If we keep increasing  the frequency our asymptotic approximation  breaks down, which physically relates   to
 the skin effect for conductive materials.
 Figure \ref{multi-recog} shows the classification results for a cube.  At each noise level, we run the algorithm 1000 times and average the results. 
This clearly illustrates that using  multiple frequencies for shape descriptors makes it possible to
distinguish a cube from a sphere.

\color{black}
\section{Concluding remarks}

In this paper we have developed an efficient classification algorithm from induction data based on dictionary matching
of shape descriptors.
This was done under the assumption that the characteristic size of the target is
of the same order of magnitude or smaller than the skin depth. The shape descriptors are constructed from conductive
 polarization tensors at multiple frequencies. If a target has a different magnetic permeability  from the
  background medium, then its second polarization tensor associated with the magnetic contrast can be extracted
  from the data and used to better classify the target. The combined use of these two polarization tensors for
  classification will be the subject of a forthcoming publication.
In future work, we will also investigate the effect of medium noise
on the classification capabilities of our proposed multifrequency, induction based, algorithm.
Our  algorithm is currently limited to the case of 
 well separated targets. 
Extending it to the case of clustered objects will likely prove to be quite challenging.
Since that case is of great importance in some practical applications, we will certainly study it 
at some point in the future.

\end{document}